\documentclass{cedram-afst}

\usepackage{bm}
\usepackage[T1]{fontenc}
\usepackage[english]{babel}
\usepackage[utf8]{inputenc}
\usepackage[mathscr]{euscript}
\usepackage{textcomp,multicol,enumitem}
\usepackage[papersize={195mm,250mm},top=2cm, bottom=2cm, left=2cm , right=2cm]{geometry}

% ...

% ...

\title
[Algebraic independence]
{On the Algebraic Independence of $E$- and $G$-Functions, I: A $p$-adic Criterion}

\alttitle{}

\author{\firstname{Daniel} \middlename{} \lastname{vargas-Montoya}}

\urladdr{}

\thanks{}

  
\altkeywords{}


\begin{document}
%% Abstracts must be placed before \maketitle
\begin{abstract}
%In this work we provide a new method for studying the algebraic independence for a large class of $G$\nobreakdash-functions. This method is based on the fact that for a lot of $G$-function $f(z)$, the logarithmic derivative of $f(z)$ belongs to $E_p$. So, thanks to a result due to Kolchin, if $f %
%Let \( f_1(z), \ldots, f_m(z) \) be power series in \( \mathbb{Q}_p[[z]] \) such that, for every \( 1 \leq i \leq m \), \( f_i(z) \) is a solution of a differential operator \( \mathcal{L}_i \in E_p[d/dz] \), where \( E_p \) is the field of analytic elements. We prove that if, for every \( 1 \leq i \leq m \), \( \mathcal{L}_i \) has a strong Frobenius structure and has maximal order multiplicity at zero (MOM), then \( f_1(z), \ldots, f_m(z) \) are algebraically dependent over \( E_p \) if and only if there are integers \( a_1, \ldots, a_m \), not all zero, such that \( f_1^{a_1}(z) \cdots f_m^{a_m}(z) \in E_p \). The main consequence of this result is that it allows us to study the algebraic independence of a large class of \( G \)-functions and certain \( E \)-functions.

Let $K$ be a finite extension of $\mathbb{Q}_p$, and let $f_1(z),\ldots, f_m(z) \in K[[z]]$ be power series such that, for every $1 \leq i \leq m$, $f_i(z)$ is a solution to a differential operator $\mathcal{L}_i \in E_p[\mathrm{d}/\mathrm{d}z]$, where $E_p$ denotes the field of analytic elements. Suppose that $K$ is totally ramified over $\mathbb{Q}_p$, and that for every $1 \leq i \leq m$, the operator $\mathcal{L}_i$ has a strong Frobenius structure and satisfies the maximal order multiplicity (MOM) condition at zero. Under these assumptions, we show that $f_1(z),\ldots, f_m(z)$ are algebraically dependent over $E_p$ if and only if there exist integers $a_1,\ldots, a_m$, not all zero, such that $f_1(z)^{a_1} \cdots f_m(z)^{a_m} \in E_p.$ The main consequence of this result is that it provides a general criterion for studying the algebraic independence of a broad class of $G$-functions and certain $E$-functions over the field of analytic elements.

\end{abstract}

%% French abstract
%\begin{altabstract}
%Ceci est le r\'esum\'e fran\c cais.
%\end{altabstract}

\maketitle

%\tableofcontents

\section{Introduction}

The classes of $E$- and $G$-functions were introduced by Siegel in his seminal 1929 paper~\cite{Siegel}. Siegel's primary aim in introducing these classes was to generalize the classical theorems of Hermite and Lindemann--Weierstrass. Recall that a power series $f(z)=\sum_{n\geq0}a_nz^n$ is a $G$-\emph{function} if the coefficients $a_n$ are algebraic numbers and there exists a real number $C>0$ such that:

1. the power series $f(z)$ is a solution to a nonzero differential operator with coefficients in $\overline{\mathbb{Q}}(z)$;\smallskip

2. the absolute values of all Galois conjugates of $a_n$ are bounded by $C^{n+1}$ for all $n\geq0$;\smallskip

3. there exists a sequence of positive integers $D_m$ such that, for all integers $m\geq0$, $D_m < C^{m+1}$ and, for all $n\leq m$, $D_m a_n$ is an algebraic integer.

Examples of $G$-functions include the hypergeometric series $_{n}F_{n-1}$ with rational parameters, as well as the diagonals of rational functions.

Furthermore, a power series $F(z)=\sum_{n\geq0}\frac{a_n}{n!}z^n$ is an $E$-\emph{function} if and only if $f(z)=\sum_{n\geq0}a_nz^n$ is a $G$-function. In particular, this implies that if $F(z)$ is an $E$-function, then it satisfies a nonzero linear differential equation with coefficients in $\overline{\mathbb{Q}}(z)$. Classical examples of $E$-functions include the exponential, sine, and cosine functions, as well as the hypergeometric series $_{p}F_p$ with rational parameters and Bessel functions.

The study of the algebraic independence of $E$- and $G$-functions is deeply motivated by Siegel's work~\cite{Siegel}, in which he developed a method to prove the algebraic independence of values of $E$\nobreakdash-functions at algebraic points. A major milestone in this theory is given by the \emph{Siegel--Shidlovskii Theorem}. Thanks to this theorem, studying the algebraic independence of $E$-functions at algebraic points is essentially equivalent to studying the algebraic independence of the functions themselves. Siegel also suggested that his method could be applied to $G$-functions. Indeed, under certain conditions, Bombieri~\cite{B} and Galochkin~\cite{G} successfully adapted Siegel's approach to obtain results on the values of $G$-functions at algebraic points. In 1984, Chudnovsky~\cite{Chu} succeeded in removing these extra conditions. However, these results concern \emph{irrationality} rather than \emph{transcendence}. Furthermore, it follows from \cite{BW} that the Siegel--Shidlovskii Theorem does not hold for $G$-functions. Nevertheless, over the past few years, several results on transcendence (see, e.g., \cite{allouche, BH, transcedencia}) and algebraic independence of $G$-functions have been established (see, e.g., \cite{ABD19,vargas2}). In contrast, the transcendence and algebraic independence of the values of $G$-functions at algebraic points remain largely open.\footnote{For more details on what is known regarding the transcendence of values of $G$-functions at algebraic points, we refer the reader to \cite[Section 7]{Ri}.}

Since $E$- and $G$-functions are solutions to differential operators with coefficients in $\overline{\mathbb{Q}}(z)$, \emph{differential Galois theory} provides a natural framework for studying their algebraic independence. For instance, using this theory, certain criteria for the algebraic independence of $E$-functions were established in \cite{B88} and \cite{shi}. Furthermore, Beukers and Heckman~\cite{BH} employed this approach to characterize the hypergeometric series ${}_{n}F_{n-1}$ that are algebraic over $\mathbb{C}(z)$. To study the algebraic independence of $E$\nobreakdash- and $G$\nobreakdash-functions $f_1(z), \ldots, f_m(z)$, one typically seeks to compute the differential Galois group associated with a differential operator annihilating these series. In general, determining this group is quite challenging. Nevertheless, a simpler case arises when each power series is a solution to a first\nobreakdash-order differential operator. In such cases, Kolchin’s Theorem provides a complete description of the algebraic relations among them. 

Our key observation is that there exist $E$\nobreakdash- and $G$-functions satisfying a first-order differential operator not over $\mathbb{Q}(z)$, but over the field of \emph{analytic elements} $E_p$. More precisely, we introduce the set $\mathcal{M}\mathcal{F}(K)$, where $K$ is a totally ramified finite extension of $\mathbb{Q}_p$. This set consists of power series $f(z)\in 1 + zK[[z]]$ such that $f(z)$ is a solution to a differential operator $\mathcal{L}\in E_p[\mathrm{d}/\mathrm{d}z]$ having \emph{maximal order multiplicity} (MOM) at zero and admitting a \emph{strong Frobenius structure}. We then establish the following result.

\begin{theo}[Theorem~\ref{theo_analytic_element}]\label{theo_weak_analytic}
Let $K$ be a finite extension of $\mathbb{Q}_p$ and let $f(z) \in 1+zK[[z]]$. If $K$ is totally ramified over $\mathbb{Q}_p$ and $f(z)\in\mathcal{M}\mathcal{F}(K)$, then:
\begin{enumerate}[label=(\roman*)]
\item $f(z)\in 1+z\mathcal{O}_K[[z]]$;
\item $f(z)$ is a solution to a differential operator of order $1$ with coefficients in $E_p$. In other words, $f'(z)/f(z) \in E_p$.
\end{enumerate}
\end{theo}

Furthermore, many $G$-functions and some $E$-functions belong to $\mathcal{M}\mathcal{F}(\mathbb{Q}_p(\pi_p))$\footnote{As we will see in Section~\ref{sec_sFs_MOM}, the field $\mathbb{Q}_p(\pi_p)$ is totally ramified over $\mathbb{Q}_p$.}, where $\pi_p$ satisfies $\pi_p^{p-1} = -p$. We will return to this observation in more detail in Section~\ref{sec_alg_ind_E_G}.

By combining Theorem~\ref{theo_weak_analytic} with a fundamental result from differential Galois theory due to Kolchin~\cite{kolchin} (see Theorem~\ref{theo_Kolchin}), we are able to prove the following theorem.

\begin{theo}[Theorems~\ref{theo_alg_ind} and \ref{theo_derivatives}]\label{theo_main}
Let $K$ be a finite extension of $\mathbb{Q}_p$ and let $f_1(z),\ldots, f_m(z)\in\mathcal{M}\mathcal{F}(K)$. If $K$ is totally ramified over $\mathbb{Q}_p$, then:
\begin{enumerate}[label=(\roman*)]
\item $f_1(z),\ldots, f_m(z)$ are algebraically dependent over $E_p$ if and only if there exist $a_1,\ldots, a_m\in\mathbb{Z}$, not all zero, such that 
\[
f_1(z)^{a_1}\cdots f_m(z)^{a_m}\in E_p;
\]
\item $f_1(z),\ldots, f_m(z)$ are algebraically dependent over $E_p$ if and only if, for any $(r_1,\ldots, r_m)\in\mathbb{N}^m$, $f_1^{(r_1)}(z),\ldots, f_m^{(r_m)}(z)$ are algebraically dependent over $E_p$.
\end{enumerate}
\end{theo}

Let us briefly explain why the MOM and strong Frobenius structure conditions arise naturally in the context of $G$-functions. In most cases, we are interested in $G$-functions $f(z) \in \mathbb{Q}[[z]]$ that are solutions to a differential operator $\mathcal{L} \in \mathbb{Q}(z)[\mathrm{d}/\mathrm{d}z]$.\footnote{Since $\mathbb{Q}(z) \subset E_p$ for any prime number $p$, we can always view $\mathcal{L}$ as an element of $E_p[\mathrm{d}/\mathrm{d}z]$.}
According to the Bombieri--Dwork conjecture (see André~\cite[p.~111]{andre} and Christol~\cite[Conjectures 2 and 3]{gilleglobalborne}), the differential operator $\mathcal{L}$ should come \emph{from geometry}; for this reason, it is expected that $\mathcal{L}$ has a strong Frobenius structure for almost every prime $p$. Although this conjecture remains open, many $G$-functions are known to be solutions to differential operators with a strong Frobenius structure for almost every $p$. Examples include hypergeometric series ${}_nF_{n-1}$ (see, e.g., \cite{vargas1}) and diagonals of rational functions (see \cite{picardfuchs, Kedlaya}). In addition, the MOM condition means that (i) zero is a regular singular point, and (ii) the \emph{indicial polynomial} of $\mathcal{L}$ at zero is $X^n$, where $n$ is the order of $\mathcal{L}$. Although the MOM condition is not always satisfied, many $G$-functions are solutions to differential operators that satisfy MOM at zero. As a primary class of $G$-functions satisfying both conditions, we have those generated by Apéry\nobreakdash-like sequences (see \cite{apery} for a comprehensive list of such sequences). Moreover, among the more than 400 operators listed in \cite{calabi}, many satisfy both conditions, and their solutions are $G$-functions.

Another interesting observation is that there exist $E$-functions that are also solutions to differential operators admitting a strong Frobenius structure. In contrast to the case of $G$-functions, as we will point out in Section~\ref{sec_alg_ind_E_G}, the family of $E$-functions satisfying this property is \emph{a priori} not as rich as that of $G$-functions. Nevertheless, series such as the exponential function $\exp(\pi_p z)$ and the Bessel function $J_0(\pi_p z)$ belong to this family. Furthermore, their corresponding differential operators also satisfy the MOM condition at zero.

%However, the MOM condition less general than the Frobenius structure. 
%for $\mathcal{L}$ is MOM means that when writing $\mathcal{L}$ in terms of $\delta=z\frac{d}{dz}$, we obtain a diffferential operator $\delta^n+a_{1}(z)\delta^{n-1}+\cdots+a_{n-1}(z)\delta+a_n(z)$ we have $a_i(z)\in\mathbb{Q}(z)\cap\mathbb{Q}[[z]]$ and $a_i(0)=0$.

%The MOM condition at zero actually means that $\mathcal{L}$ has \emph{maximal unipotent monodromy} at zero. However, since we are interested in studying 
 
%\textbf{3) In Section~\ref{sec_alg_ind_E_G}, we point out that the set $\mathcal{M}\mathcal{F}(K)$ contains an interesting class of $E$- and $G$\nobreakdash-functions. In a lot of examples, the differential operator $\mathcal{L}$ belongs to $\mathbb{Q}(z)[d/dz]$, but since $\mathbb{Q}(z)\subset E_p$, we then see $\mathcal{L}$ as an element of $E_p[d/dz]$. In addition, the MOM property at zero for $\mathcal{L}\in\mathbb{Q}(z)[d/dz]$ is equivalent to saying that $\mathcal{L}$ has maximal unipotent monodromy at zero (MUM). The notion of strong Frobenius structure was introduced by Dwork~\cite{DworksFf} and it corresponds to a Frobenius action on the set of solutions of $\mathcal{L}$. In Section~\ref{sec_alg_ind_E_G}, we give some examples of differential operator having strong Frobenius structure. }

In Part~II, we made Theorem~\ref{theo_main} effective by providing a criterion to determine when $\prod_{i=1}^m f_i(z)^{a_i} \notin E_p$ for any nonzero tuple $(a_1,\ldots, a_m) \in \mathbb{Z}^m$. This criterion, stated as Theorem~2.8 in Part~II, relies on analyzing the poles and residues of the logarithmic derivatives $f'_1(z)/f_1(z),\ldots, f'_m(z)/f_m(z)$ viewed as elements of $E_p$. In the same work, we illustrated this criterion by proving the algebraic independence of several $E$- and $G$-functions.

For instance, let us consider the power series 
\[
J_0(z)=\sum_{n\geq0}\frac{(-1)^n}{4^n(n!)^2}z^{2n}, \quad \mathfrak{F}(z)=\sum_{n\geq0}\left(\sum_{k=0}^n\binom{n}{k}^3\right)z^n, \quad \text{and} \quad \mathfrak{A}(z)=\sum_{n\geq0}\left(\sum_{k=0}^n\binom{n}{k}^2\binom{n+k}{k}^2\right)z^n.
\]
We prove in Part~II that for any non-negative integers $r, s, k \geq 0$, the power series $J_0^{(r)}(\pi_5 z)$, $\mathfrak{A}^{(s)}(z)$, and $\mathfrak{F}^{(k)}(z)$ are algebraically independent over $E_5$, where $\pi_5^2 = -5$. Since $\mathbb{Q}(z) \subset E_p$ for every prime $p$, this result immediately implies that $J_0^{(r)}(\pi_5 z)$, $\mathfrak{A}^{(s)}(z)$, and $\mathfrak{F}^{(k)}(z)$ are algebraically independent over $\mathbb{Q}(z)$. In this setup, $J_0(\pi_5 z)$ is an $E$-function, while $\mathfrak{F}(z)$ and $\mathfrak{A}(z)$ are $G$-functions.   

Furthermore, it is worth pointing out that in Part~II we also obtain results on algebraic independence over $\mathbb{C}_p(z)$. More precisely, in Theorem~2.9 of Part~II, we prove that, under certain conditions, the power series $f_1(z),\ldots, f_m(z)\in\mathcal{MF}(K)$ are algebraically independent over $E_p$ if and only if there exist integers $a_1,\ldots, a_m\in\mathbb{Z}$, not all zero, such that $f_1(z)^{a_1}\cdots f_m(z)^{a_m}\in\mathbb{C}_p(z)$.\footnote{This result is quite surprising because there exist power series that are algebraic over $E_p$ but transcendental over $\mathbb{C}_p(z)$, as is the case for the power series $\exp(\pi_p z)$.}

For example, by using this result, we can show that 
\[
\mathfrak{a}(z)=\sum_{n\geq0}\left(\sum_{k=0}^n\binom{n}{k}^2\binom{n+k}{k}\right)z^n \quad \text{and} \quad \mathfrak{A}(z)=\sum_{n\geq0}\left(\sum_{k=0}^n\binom{n}{k}^2\binom{n+k}{k}^2\right)z^n
\]
are algebraically independent over $\mathbb{C}_3(z)$.
% The proof of Theorem~\ref{theo_main} combines differential Galois theory and the theory of $p$\nobreakdash-adic differential equations. Using the latter, we prove the following result:

 %We refer the reader to Section~\ref{sec_proof_results} for a precise statement of Kolchin's Theorem and how we apply it in order to prove Theorem~\ref{theo_main}.
  
%As an illustration of Theorem~\ref{theo_main} we prove the following result
% \begin{theo}[Theorem \ref{theo_apery_bessel}]
%Let us consider the Bessel fonction $$J_0(z)=\sum_{n\geq0}\frac{(-1)^n}{4^n(n!)^2}z^{2n}.$$
%Then, for any prime number $p>2$ and any integer $r\geq0$, the power series $J_0^{(r)}(\pi_pz)$ is transcendental over $E_3$.
%\end{theo}
%It is worth mentioning that this last result confirms a conjecture due to Christol (see \cite[p.30]{C86}).
% This result implies that,  for all integers $r,s, k\geq0$, the power series $J_0^{(r)}(\pi_3z)$, $\mathfrak{A}^{(s)}(z)$, and $\mathfrak{f}^{(k)}(z)$ are algebraically independent over $\mathbb{Q}(z)$ because, for all prime numbers $p$, $\mathbb{Q}(z)\subset E_p$. In this case, the power series $J_0(\pi_3z)$ is a $E$-function and the power series $\mathfrak{f}(z)$ and $\mathfrak{A}(z)$ are $G$-functions.

A result similar to part~(i) of Theorem~\ref{theo_main} and Theorem~2.9 of Part~II was proven by Adamczewski, Bell, and Delaygue~\cite{ABD19} for the set of power series that are $p$-Lucas for infinitely many prime numbers $p$. Recall that a power series $f(z)=\sum_{n\geq0}a_nz^n\in\mathbb{Q}[[z]]$ is \emph{$p$-Lucas} if $a_0=1$, $f(z) \in \mathbb{Z}_{(p)}[[z]]$\footnote{Recall that $\mathbb{Z}_{(p)}$ is the ring of rational numbers $a/b$ such that $p$ does not divide $b$ and $a$ is coprime to $b$.}, and 
\begin{equation*}
f(z) \equiv \left(\sum_{n=0}^{p-1}a_nz^n\right)f(z^p) \bmod{p\mathbb{Z}_{(p)}[[z]]}.
\end{equation*}
In \cite{ABD19}, the authors introduced the set of power series $\mathcal{L}(\mathcal{S})$, where $\mathcal{S}$ is an infinite set of prime numbers. The main examples in this set consist of power series that are $p$-Lucas for all primes $p \in \mathcal{S}$. The main result of \cite{ABD19} states that if $f_1(z),\ldots, f_m(z)$ belong to $\mathcal{L}(\mathcal{S})$, then $f_1(z),\ldots, f_m(z)$ are algebraically independent over $\mathbb{Q}(z)$ if and only if there exist $a_1,\ldots, a_m\in\mathbb{Z}$, not all zero, such that $f_1(z)^{a_1}\cdots f_m(z)^{a_m}\in\mathbb{Q}(z)$. In the approach developed in~\cite{ABD19}, assuming that $\mathcal{S}$ is infinite is crucial. In the same work, the authors also showed that a large class of $G$-functions are $p$-Lucas for infinitely many primes $p$. Furthermore, in \cite{vargas2}, we introduced the set $\mathcal{L}^2(\mathcal{S})$. This set contains $\mathcal{L}(\mathcal{S})$, and we exhibited $G$-functions that belong to $\mathcal{L}^2(\mathcal{S})$ but not to $\mathcal{L}(\mathcal{S})$ for any infinite set $\mathcal{S}$ of prime numbers. Moreover, a criterion for the algebraic independence of power series in $\mathcal{L}^2(\mathcal{S})$ was established and illustrated by proving the algebraic independence of several $G$-functions that are not $p$-Lucas for any prime number $p$ (see, e.g., \cite[Theorems 9.1, 9.2]{vargas2}).

It is worth mentioning that, based on numerous examples, we observe that $G$-functions belonging to $\mathcal{L}(\mathcal{S})$ or $\mathcal{L}^2(\mathcal{S})$ also lie in $\mathcal{M}\mathcal{F}(\mathbb{Q}_p)$ for almost every prime number $p \in \mathcal{S}$. Another important point is that if $f(z) \in \mathcal{L}(\mathcal{S})$, then typically $f'(z) \notin \mathcal{L}(\mathcal{S})$. Nevertheless, the strategy used in \cite{ABD19} seems adaptable to studying the algebraic independence of derivatives of elements in $\mathcal{L}(\mathcal{S})$. In contrast, the method developed in this work provides an alternative approach to studying the algebraic independence of a large class of $G$-functions and their derivatives. A key advantage of our approach is that it focuses on a single prime number $p$. This feature makes it possible to include $E$-functions in our study, something that \emph{a priori} is not possible with the results of \cite{ABD19, vargas2}, since the sets $\mathcal{L}(\mathcal{S})$ and $\mathcal{L}^2(\mathcal{S})$ do not contain any $E$-functions. Moreover, our method also allows us to study simultaneously the algebraic independence of $E$- and $G$-functions.

The outline of the paper is as follows. In Section~\ref{sec_results}, we recall the definitions of the field of analytic elements and the notion of strong Frobenius structure. In Section~\ref{sec_proof_results}, we prove Theorem~\ref{theo_main}. In Section~\ref{sec_stra}, we present the strategy for proving Theorem~\ref{theo_weak_analytic}, highlighting the main steps of the proof. Finally, Sections~\ref{sec_frob_ant} and \ref{sec_congruences} are devoted to developing the strategy outlined in Section~\ref{sec_stra}. %Finally, in Section~\ref{sec_appli}  we provide a necessary condition for determining the  algebraic independence over $E_{p}$. This condition is given by Theorem~\ref{theo_criterion_alg_ind} and we illustrate this criterion by proving the algebraic independence of some $G$\nobreakdash-functions and $E$\nobreakdash-functions. 

%\textbf{Acknowledgements}. We warmly thank Charlotte Hardouin for the valuable discussions regarding this project. The second part of Theorem~\ref{theo_main} was motivated by a question from Tanguy Rivoal about whether our method could be used to study the algebraic independence of derivatives; we also thank him.

 %Moreover, in \cite{vargas1} it is proven that $\sum_{j=0}^{n^2h}a_i(z)f_{\mid p}(z)^{p^{j}}=0.$
 %where $a_j(z)\in\mathbb{F}_p(z)$, $n$ is the order of $L$, and $h$ is the period of the strong Frobenius structure. Notice that if $g(z)$ belongs to $\mathcal{L}(\mathcal{S})$ then, thanks to \eqref{eq_p_lucas_g}, $g_{\mid p}(z)$ is algebraic over $\mathbb{F}_p(z)$

 \section{Main results}\label{sec_results}
 
We first recall the definitions of the field of analytic elements, strong Frobenius structures, and MOM operators.

\subsection{Analytic elements}
Let $p$ be a prime number, and let $\mathbb{Q}_p$ be the field of $p$-adic numbers. It is well known that $\mathbb{Q}_p$ is equipped with the non-Archimedean norm $|\cdot|$ such that $|p|=1/p$; this norm is called the $p$-adic norm. Let $\overline{\mathbb{Q}_p}$ be an algebraic closure of $\mathbb{Q}_p$. This norm extends uniquely to $\overline{\mathbb{Q}_p}$ and is still called the $p$-adic norm. Finally, let $\mathbb{C}_p$ denote the completion of $\overline{\mathbb{Q}_p}$ with respect to the $p$-adic norm.

From now on, $K \subset \mathbb{C}_p$ will denote a complete extension of $\mathbb{Q}_p$ with respect to the $p$-adic norm.

Recall that the Amice ring with coefficients in $K$ is defined as follows: 
\[
\mathcal{A}_K = \left\{\sum_{n\in\mathbb{Z}}a_n z^n : a_n \in K, \lim\limits_{n\to -\infty}|a_n|=0 \text{ and } \sup\limits_{n\in\mathbb{Z}}|a_n|<\infty\right\}.
\]
This ring is equipped with the Gauss norm, defined by: 
\[
\left|\sum_{n\in\mathbb{Z}}a_n z^n\right|_{\mathcal{G}} = \sup\limits_{n\in\mathbb{Z}}|a_n|.
\]
Since $K$ is complete, Proposition~1.1 of \cite{C86} implies that the ring $\mathcal{A}_K$ is also complete with respect to the Gauss norm. Furthermore, it follows from Proposition~1.2 of \cite{C86} that an element $h = \sum_{n\in\mathbb{Z}} h_n z^n \in \mathcal{A}_K$ is a unit if and only if there exists $n_0 \in \mathbb{Z}$ such that $|h|_{\mathcal{G}} = |h_{n_0}| \neq 0$. In particular, any nonzero polynomial with coefficients in $K$ is a unit in $\mathcal{A}_K$, and therefore $K(z) \subset \mathcal{A}_K$. For this reason, the field $K(z)$ is equipped with the Gauss norm. In fact, for any nonzero polynomials $P(z) = \sum_{i=0}^n a_i z^i$ and $Q(z) = \sum_{j=0}^m b_j z^j$, we have 
\[
\left|\frac{P(z)}{Q(z)}\right|_{\mathcal{G}} = \frac{\max_{0 \leq i \leq n} |a_i|}{\max_{0 \leq j \leq m} |b_j|}.
\]
The field of \emph{analytic elements} with coefficients in $K$, denoted by $E_K$, is defined as the completion of $K(z)$ with respect to the Gauss norm. Since $\mathcal{A}_K$ is complete and $K(z) \subset \mathcal{A}_K$, we have $E_K \subset \mathcal{A}_K$. When $K = \mathbb{C}_p$, the field $E_{\mathbb{C}_p}$ is denoted by $E_p$, and the Amice ring $\mathcal{A}_{\mathbb{C}_p}$ is denoted by $\mathcal{A}_p$.

\begin{rema}\label{rem_contenido}
We have $\mathbb{Q}(z) \subset \mathbb{Q}_p(z) \subset K(z) \subset E_K$ because $\mathbb{Q} \subset \mathbb{Q}_p \subset K$, and clearly $K(z) \subset E_K$. Similarly, $E_K \subset E_{K'}$ if $K \subset K'$.
\end{rema}

Finally, consider the open unit disk $D_0 = \{x \in \mathbb{C}_p : |x| < 1\}$. We denote by $K_0(z)$ the ring of rational functions $A(z)/B(z) \in K(z)$ such that $B(x) \neq 0$ for all $x \in D_0$. Let $E_{0,K}$ be the completion of $K_0(z)$ with respect to the Gauss norm. It is clear that $E_{0,K} \subset E_K$ and $E_{0,K} \subset K[[z]]$. When $K = \mathbb{C}_p$, the ring $E_{0,\mathbb{C}_p}$ is denoted by $E_{0,p}$.

\subsection{Strong Frobenius structures and MOM operators}\label{sec_sFs_MOM}

We recall that two matrices $B_1, B_2 \in M_n(E_K)$ are \emph{$E_K$-equivalent} if there exists $T \in \mathrm{GL}_n(E_p)$ such that $\delta T = B_1 T - T B_2$, where $\delta = z\frac{\mathrm{d}}{\mathrm{d}z}$. We also recall that an endomorphism $\sigma \colon K \to K$ is a \emph{Frobenius endomorphism} if $\sigma$ is continuous and, for all $x \in K$ with $|x| \leq 1$, we have $|\sigma(x) - x^p| < 1$. 

Let us assume that $K$ is equipped with a Frobenius endomorphism $\sigma$. We define the map $\bm\sigma \colon K(z) \to E_K$ by 
\[
\bm\sigma \left( \frac{\sum_{i=0}^l a_i z^i}{\sum_{j=0}^s b_j z^j} \right) = \frac{\sum_{i=0}^l \sigma(a_i) z^{i p}}{\sum_{j=0}^s \sigma(b_j) z^{j p}}.
\]
Since $\sigma$ is continuous, $\bm\sigma$ is continuous. Therefore, $\bm\sigma$ extends uniquely to a continuous endomorphism of $E_K$, which we also denote by $\bm\sigma$. 

Let $L = \delta^n + a_1(z)\delta^{n-1} + \cdots + a_{n-1}(z)\delta + a_n(z)$ be a differential operator in $E_K[\delta]$. The companion matrix of $L$ is 
\[
A = \begin{pmatrix}
0 & 1 & 0 & \dots & 0 & 0\\
0 & 0 & 1 & \dots & 0 & 0\\
\vdots & \vdots & \vdots & \ddots & \vdots & \vdots \\
0 & 0 & 0 & \dots & 0 & 1\\
-a_n(z) & -a_{n-1}(z) & -a_{n-2}(z) & \dots & -a_2(z) & -a_1(z)
\end{pmatrix}.
\]
We recall that $L$ has a \emph{strong Frobenius structure} of period $h$ if the matrices $A$ and $p^h \bm\sigma^h(A)$ are $E_p$-equivalent, where $\bm\sigma^h(A)$ denotes the matrix obtained by applying $\bm\sigma^h$ entrywise to $A$. Notice that if $K$ is totally ramified over $\mathbb{Q}_p$ (that is, its residue field is $\mathbb{F}_p$), then the identity map $\mathrm{id} \colon K \to K$ is a Frobenius endomorphism. From now on, we say that $K$ is a \emph{Frobenius field} if $K$ is a finite extension of $\mathbb{Q}_p$ and the identity map $\mathrm{id} \colon K \to K$ is a Frobenius endomorphism. Therefore, if $K$ is a Frobenius field, then $L \in E_K[\delta]$ has a strong Frobenius structure of period $h$ if the matrices $A$ and $p^h A(z^{p^h})$ are $E_p$-equivalent. Here are two important examples of Frobenius fields:

\begin{itemize}
\item The field $\mathbb{Q}_p$ is a Frobenius field because its residue field is $\mathbb{F}_p$.
\item Any root of $X^{p-1} + p$ is called a \emph{Dwork constant} and is denoted by $\pi_p$. The field $\mathbb{Q}_p(\pi_p)$ is a Frobenius field. Indeed, since $\pi_p$ is a root of $X^{p-1} + p$ and this polynomial is Eisenstein, it follows that $\mathbb{Q}_p(\pi_p)$ is totally ramified over $\mathbb{Q}_p$, and therefore $\mathbb{Q}_p(\pi_p)$ is a Frobenius field.
\end{itemize}

Let $L \in \mathbb{Q}(z)[\delta]$, and let $A$ be its companion matrix. Remark~\ref{rem_contenido} implies that $L \in E_{\mathbb{Q}_p}[\delta]$. We then say that $L$ has a \emph{strong Frobenius structure} of period $h$ for $p$ if, viewed as an element of $E_{\mathbb{Q}_p}[\delta]$, $L$ has a strong Frobenius structure of period $h$. Since $\mathbb{Q}_p$ is a Frobenius field, this is equivalent to saying that $A$ and $p^h A(z^{p^h})$ are $E_p$-equivalent.

Finally, we say that $L = \delta^n + a_1(z)\delta^{n-1} + \cdots + a_{n-1}(z)\delta + a_n(z) \in E_K[\delta]$ is MOM (\emph{maximal order of multiplicity}) at zero if, for all $i \in \{1,\ldots, n\}$, we have $a_i(z) \in K[[z]]$ and $a_i(0) = 0$. In fact, when $L \in \mathbb{Q}(z)[\delta]$, being MOM at zero is equivalent to saying that $L$ is MUM (\emph{maximal unipotent monodromy}) at zero.

 %According to Remark~\ref{rem_sol}, if $L$ is MOM at zero then $L$ has a solution in $1+zK[[z]]$. 

\subsection{Main results}

We first introduce the set of power series $\mathcal{M}\mathcal{F}(K)$.

\begin{defi}
Let $K \subset \mathbb{C}_p$ be a complete extension of $\mathbb{Q}_p$. We denote by $\mathcal{M}\mathcal{F}(K)$ the set of power series $f(z) \in 1 + zK[[z]]$ such that $f(z)$ is a solution to a monic differential operator 
\[
L = \delta^n + a_1(z)\delta^{n-1} + \cdots + a_{n-1}(z)\delta + a_n(z) \in E_{0,K}[\delta]
\]
satisfying the following properties:
\begin{enumerate}
\item $L$ is MOM at zero;
\item $L$ has a strong Frobenius structure;
\item $|a_i(z)|_{\mathcal{G}} \leq 1$ for all $1 \leq i \leq n$.
\end{enumerate}
\end{defi}

We are now ready to state our main results.

\begin{theo}\label{theo_alg_ind}
Let $K$ be a Frobenius field, and let $f_1(z),\ldots, f_m(z) \in \mathcal{M}\mathcal{F}(K)$. Then $f_1(z),\ldots, f_m(z)$ are algebraically dependent over $E_K$ if and only if there exist integers $a_1,\ldots, a_m \in \mathbb{Z}$, not all zero, such that 
\[
f_1(z)^{a_1}\cdots f_m(z)^{a_m} \in E_{0,K}.
\]
\end{theo}

The main ingredients in the proof of Theorem~\ref{theo_alg_ind} are Theorem~\ref{theo_weak_analytic} and Kolchin's Theorem (see Theorem~\ref{theo_Kolchin}). Although an element $f(z) \in \mathcal{M}\mathcal{F}(K)$ does not necessarily satisfy $f'(z) \in \mathcal{M}\mathcal{F}(K)$, the strategy used to prove Theorem~\ref{theo_alg_ind} also allows us to study the algebraic independence of its higher derivatives. More precisely:

\begin{theo}\label{theo_derivatives}
Let $K$ be a Frobenius field, and let $f_1(z),\ldots, f_m(z) \in \mathcal{M}\mathcal{F}(K)$ be such that, for all $1 \leq i \leq m$, $f_i^{(r)}(z) \neq 0$ for every integer $r \geq 0$. Then the following statements are equivalent:
\begin{enumerate}[label=(\roman*)]
\item $f_1(z),\ldots, f_m(z)$ are algebraically dependent over $E_K$;
\item for all $(r_1,\ldots, r_m) \in \mathbb{N}^m$, the derivatives $f_1^{(r_1)}(z),\ldots, f_m^{(r_m)}(z)$ are algebraically dependent over $E_K$;
\item for some $(s_1,\ldots, s_m) \in \mathbb{N}^m$, the derivatives $f_1^{(s_1)}(z),\ldots, f_m^{(s_m)}(z)$ are algebraically dependent over $E_K$.
\end{enumerate}
\end{theo}

\subsection{Algebraic independence of $G$-functions and $E$-functions over $E_p$}\label{sec_alg_ind_E_G}

Theorems~\ref{theo_alg_ind} and \ref{theo_derivatives} put us in a position to study the algebraic independence of $G$\nobreakdash-functions and $E$\nobreakdash-functions over the field $E_p$, since many $G$\nobreakdash-functions belong to $\mathcal{M}\mathcal{F}(\mathbb{Q}_p)$ and certain $E$\nobreakdash-functions belong to $\mathcal{M}\mathcal{F}(\mathbb{Q}_p(\pi_p))$. Let us examine this in more detail.

\textbullet\quad\textbf{$G$-functions}

The most representative examples of $G$-functions are the \emph{hypergeometric series} ${}_nF_{n-1}$ and the \emph{diagonals of rational functions}. In what follows, we prove that many of these power series belong to $\mathcal{M}\mathcal{F}(\mathbb{Q}_p)$. 

\paragraph{Hypergeometric series ${}_nF_{n-1}$.} 
Recall that the hypergeometric series with parameters $\bm{\alpha}=(\alpha_1,\ldots, \alpha_n)$ and $\bm{\beta}=(\beta_1,\ldots, \beta_{n-1},1)$ in $(\mathbb{Q}\setminus\mathbb{Z}_{\leq0})^n$ is given by 
\[
{}_nF_{n-1}(\bm{\alpha},\bm{\beta},z) = \sum_{j\geq0}\frac{(\alpha_1)_j\cdots(\alpha_n)_j}{(\beta_1)_j\cdots(\beta_{n-1})_j\, j!}z^j,
\]
where $(x)_j$ denotes the Pochhammer symbol, defined as $(x)_0=1$ and $(x)_j=x(x+1)\cdots(x+j-1)$ for $j>0$. Furthermore, it is well known that ${}_nF_{n-1}(\bm{\alpha},\bm{\beta},z)$ is a solution to the differential operator 
\[
\mathcal{H}(\bm{\alpha},\bm{\beta}) = \prod_{i=1}^n(\delta+\beta_i-1) - z\prod_{i=1}^n(\delta+\alpha_i).
\]
It was proven in \cite{vargas1} that if $\alpha_i-\beta_j \notin \mathbb{Z}$ for all $1 \leq i,j \leq n$, then $\mathcal{H}(\bm{\alpha},\bm{\beta})$ has a strong Frobenius structure for all $p > d_{\bm{\alpha},\bm{\beta}}$, where $d_{\bm{\alpha},\bm{\beta}}$ is the least common multiple of the denominators of the $\alpha_i$'s and $\beta_j$'s. Moreover, if $\bm{\beta}=(1,\ldots, 1)$, then $\mathcal{H}(\bm{\alpha},\bm{\beta})$ is MOM at zero. Consequently, hypergeometric series of the form
\[
f_{\bm{\alpha}} = \sum_{j\geq0}\frac{(\alpha_1)_j\cdots(\alpha_n)_j}{(j!)^n} z^j
\]
belong to $\mathcal{M}\mathcal{F}(\mathbb{Q}_p)$ for all $p > d_{\bm{\alpha},\bm{\beta}}$, provided that $\alpha_i-1 \notin \mathbb{Z}$ for all $1 \leq i \leq n$. For example, 
\[
{}_2F_1((-1/2,1/2),(1,1),z) = \sum_{n\geq0}\frac{-1}{(2n-1)64^n}\binom{2n}{n}^2 z^n
\]
belongs to $\mathcal{M}\mathcal{F}(\mathbb{Q}_p)$ for all $p > 2$. Moreover, this power series cannot be studied using the techniques developed in \cite{ABD19}, since it was proved in \cite{vargas2} that it does not belong to $\mathcal{L}(\mathcal{S})$ for any infinite set $\mathcal{S}$ of prime numbers.

\paragraph{Diagonals of rational functions.} 
A power series $f(z) = \sum_{n\geq0} a_n z^n \in \mathbb{Q}[[z]]$ is a diagonal of a rational function if there exist an integer $m > 0$ and a rational function
\[
F(z_1,\ldots,z_m) = \sum_{(i_1,\ldots, i_m)\in\mathbb{N}^m} a_{(i_1,\ldots,i_m)} z_1^{i_1}\cdots z_m^{i_m} \in \mathbb{Q}[[z_1,\ldots, z_m]] \cap \mathbb{Q}(z_1,\ldots, z_m)
\]
such that $a_{(n,\ldots, n)} = a_n$ for all integers $n \geq 0$. Classic examples of diagonals of rational functions include:
\[
\mathfrak{A}(z) = \sum_{n\geq0}\left(\sum_{k=0}^n \binom{n}{k}^2 \binom{n+k}{k}^2 \right) z^n, \quad 
\mathfrak{A}_1(z) = \sum_{n\geq0}\left(\sum_{k=0}^n \binom{n}{k}^2 \binom{2k}{n}^2 \right) z^n,
\]
\[
\mathfrak{A}_2(z) = \sum_{n\geq0}\left(\sum_{k=0}^n \binom{n}{k}^2 \binom{2k}{k}^2 \binom{2n-2k}{n-k} \right) z^n.
\]
These $G$-functions are generated by Apéry-type sequences; a complete list of such sequences can be found in \cite{apery}. Moreover, \cite[Theorem~3.5]{sb} implies that any power series generated by an Apéry-type sequence is the diagonal of a rational function.

Another source of examples of diagonals of rational functions is the list given in \cite{calabi}. More than 400 power series are listed there, and thanks to \cite[Theorem~3.5]{sb}, many of these series are diagonals of rational functions. As an illustration, consider series \#210:
\[
\sum_{n\geq0}\binom{2n}{n}\left(\sum_{k=0}^{2n}(-1)^k\binom{2n}{k}^4\right)z^n.
\]

An important observation is that each of the $G$-functions mentioned above is a solution to a differential operator that is MOM at zero and admits a strong Frobenius structure for almost every prime number $p$. For example, $\mathfrak{A}(z)$, $\mathfrak{A}_1(z)$, and $\mathfrak{A}_2(z)$ are solutions, respectively, of the operators
\[
\delta^3 - z(2\delta+1)(17^2\delta^2+17\delta+5) + z^2(\delta+1)^3,
\]
\[
\delta^3 - z(2\delta+1)(12^2\delta^2+12\delta+4) + 16z^2(\delta+1)^3,
\]
\[
\delta^3 - z(2\delta+1)(10^2\delta^2+10\delta+4) + 64z^2(\delta+1)^3,
\]
which satisfy both properties.

The following remark provides a generalization of the previous situation.

\begin{rema}\label{rem_impo}
Suppose that $f(z) \in \mathbb{Q}[[z]]$ is a diagonal of a rational function and that $\mathcal{D}(f) = 0$ for some $\mathcal{D} \in \mathbb{Q}(z)[\delta]$. If $\mathcal{D}$ is MOM at zero, then $f(z)$ belongs to $\mathcal{M}\mathcal{F}(\mathbb{Q}_p)$ for almost every prime $p$.
\end{rema}

To prove this remark, we first recall the notion of a \emph{Picard--Fuchs} equation. Then it is well-known that, if $X$ is a smooth $\mathbb{Q}(z)$-variety, then the de Rham cohomology groups $H^i_{\mathrm{dR}}(X)$ are finite-dimensional $\mathbb{Q}(z)$-vector spaces equipped with a canonical connection $\nabla$. Therefore, for any $\omega \in H^i_{\mathrm{dR}}(X)$, there exist polynomials $a_0(z), \ldots, a_n(z) \in \mathbb{Q}[z]$, not all zero, such that 
\[
a_0(z)\omega + a_1(z)\nabla(\omega) + \cdots + a_n(z)\nabla^n(\omega) = 0.
\]
This gives rise to the differential operator $\mathcal{L}_{\omega} = a_0(z) + a_1(z)\delta + \cdots + a_n(z)\delta^n$. Following \cite[Chap.~II, Sec.~1.2]{andre}, we say that $\mathcal{L} \in \mathbb{Q}(z)[\delta]$ is a Picard--Fuchs equation if there exist a smooth $\mathbb{Q}(z)$-variety $X$ and an element $\omega$ in one of its de Rham cohomology groups such that $\mathcal{L} = \mathcal{L}_\omega$. Moreover, it follows from \cite[Theorem~22.2.1]{Kedlaya} that any Picard--Fuchs equation is equipped with a strong Frobenius structure for almost every prime $p$.

\begin{proof}[Proof of Remark~\ref{rem_impo}]
Let $\mathcal{R} \in \mathbb{Q}(z)[\delta]$ be the minimal differential operator of $f(z)$, and let $N$ be the differential module over $E_p$ associated with $\mathcal{R}$. Our goal is to prove that $\mathcal{R}$ is MOM at zero and has a strong Frobenius structure for almost every prime number $p$.

By the minimality of $\mathcal{R}$, we have $\mathcal{D} \in \mathbb{Q}(z)[\delta]\mathcal{R}$. Since $\mathcal{D}$ is MOM at zero, $\mathcal{R}$ is also MOM at zero. We now show that $N$ has a strong Frobenius structure for almost every prime $p$.

Since $f(z)$ is the diagonal of a rational function, a result of Christol~\cite{picardfuchs} implies that $f(z)$ is a solution to a Picard--Fuchs equation $\mathcal{L} \in \mathbb{Q}(z)[\delta]$; thus, $\mathcal{L}$ is equipped with a strong Frobenius structure for almost every $p$. Let $M$ be the differential module over $E_p$ associated with $\mathcal{L}$. According to \cite[Theorem~4.2.6]{delignehodge} or \cite[Proposition~7.1]{C86}, $M$ is semisimple. Since $N$ is a submodule of $M$, $N$ is also semisimple. In particular, $N = \bigoplus_{i=1}^r N_i$, where the $N_i$'s are simple modules. 

Furthermore, $M$ has a strong Frobenius structure for almost every prime $p$. Let $p$ be one of these prime numbers, and denote by $\phi_p$ the Frobenius functor, so that $\phi_p^h(M) \cong M$, where $h$ is the period of the Frobenius structure. By Krull--Schmidt's Theorem (see, e.g., \cite[p.~83]{CR66}), there exist two distinct natural numbers $m$ and $m'$ such that $\phi_p^m(N_i)$ and $\phi_p^{m'}(N_i)$ are isomorphic as differential modules over $E_p$. Hence, $N_i$ has a strong Frobenius structure for $p$. Finally, $N$ also has a strong Frobenius structure for $p$ because $\phi_p$ respects direct sums.
\end{proof}

  \textbullet\quad\textbf{E-functions}
  
\subsubsection*{$E$-functions}

We now give some examples of $E$\nobreakdash-functions belonging to $\mathcal{M}\mathcal{F}(\mathbb{Q}_p(\pi_p))$, where $\pi_p$ is Dwork's constant. Our first example is 
\[
\exp(\pi_p z) = 1 + \frac{\pi_p}{1!}z + \frac{\pi_p^2}{2!}z^2 + \cdots + \frac{\pi_p^n}{n!}z^n + \cdots.
\]
This power series is a solution to the differential operator $\delta - \pi_p z$. According to Dwork~\cite{DworksFf}, $\delta - \pi_p z$ has a strong Frobenius structure for any prime $p$. Thus, $\exp(\pi_p z)$ belongs to $\mathcal{M}\mathcal{F}(\mathbb{Q}_p(\pi_p))$. We recall that $\mathbb{Q}_p(\pi_p)$ is a Frobenius field for every prime $p$.

Another example is given by the Bessel function 
\[
J_0(\pi_p z) = \sum_{n\geq0} \frac{(-1)^n \pi_p^{2n}}{4^n (n!)^2} z^{2n}.
\]
The Bessel function $J_0(\pi_p z)$ is a solution to $\delta^2 - \pi_p^2 z^2$. It is clear that this differential operator is MOM at zero and, by Dwork~\cite{Bessel}, it admits a strong Frobenius structure for any $p \neq 2$. Thus, $J_0(\pi_p z) \in \mathcal{M}\mathcal{F}(\mathbb{Q}_p(\pi_p))$ for all $p \neq 2$. 

Further examples are provided by Dwork in \cite{DworkIV}. More precisely, \cite[Theorem~3.1]{DworkIV} implies that the logarithmic derivatives of certain hypergeometric $E$\nobreakdash-functions belong to $E_{0,K}$ with $K = \mathbb{Q}_p(\pi_p)$; consequently, these hypergeometric $E$\nobreakdash-functions lie in $\mathcal{M}\mathcal{F}(\mathbb{Q}_p(\pi_p))$.

To our knowledge, the $E$-functions mentioned above are currently the only known examples that are solutions to differential operators admitting a strong Frobenius structure. Therefore, in view of our main results, it is natural to ask under what conditions an $E$-function is a solution to a differential operator admitting a strong Frobenius structure.

\section{Proof of Theorems~\ref{theo_alg_ind} and \ref{theo_derivatives}}\label{sec_proof_results}

The proof of Theorem~\ref{theo_alg_ind} relies on Theorem~\ref{theo_analytic_element} and Kolchin's Theorem. Recall that the set $\mathcal{O}_K = \{x \in K : |x| \leq 1\}$ is a ring, commonly referred to as the ring of integers of $K$.

\begin{theo}\label{theo_analytic_element}
Let $K$ be a Frobenius field, and let $f(z) \in \mathcal{M}\mathcal{F}(K)$. Then:
\begin{enumerate}[label=(\roman*)]
\item $f(z) \in 1 + z\mathcal{O}_K[[z]]$;
\item $f(z)/f(z^{p^h}) \in E_{0,K}$, where $h$ is the period of the strong Frobenius structure of the differential operator $\mathcal{L} \in E_K[\delta]$ annihilating $f(z)$;
\item $f'(z)/f(z) \in E_{0,K}$.
\end{enumerate}
\end{theo}

It is important to point out that the assumption that $K$ is a Frobenius field is used only in the proofs of (ii) and (iii). Assertion (i), on the other hand, holds for any power series $f(z)$ that is a solution to a differential operator in $E_K[\delta]$ which is MOM at zero and admits a strong Frobenius structure. Furthermore, we have the following remark.

\begin{rema}
From the second assertion of Theorem~\ref{theo_analytic_element}, $f(z)/f(z^{p^h})$ belongs to $E_{0,K}$ and is therefore the limit of rational functions in $K_0(z)$. Nevertheless, our proof does not provide an explicit construction of rational functions converging to $f(z)/f(z^{p^h})$. This naturally raises the question: can these rational functions be determined explicitly? In certain cases, thanks to the work of Dwork~\cite{Dworkpciclos}, this is indeed possible for the hypergeometric series $f_{\bm{\alpha}}$, where $\bm{\alpha} \in (\mathbb{Q} \cap (0,1))^n$. In the particular case where $p \equiv 1 \bmod{d_{\bm{\alpha}}}$, the series $f_{\bm{\alpha}}$ satisfies the so-called $p$-Dwork congruences:
\[
\frac{f_{\bm{\alpha}}(z)}{f_{\bm{\alpha}}(z^p)} \equiv \frac{F_s(z)}{F_{s-1}(z^p)} \bmod{p^s\mathbb{Z}_{(p)}[[z]]},
\]
where $F_s$ is the $(p^s-1)$-th truncation of $f_{\bm{\alpha}}$. It turns out that the $p$-Dwork congruences have been useful in constructing the unit root of certain $p$-adic cohomologies (see, e.g., \cite{masha, BeukersMasha1, dworkfamily}). For this reason, it is of great interest to ask whether a general procedure exists to construct explicit rational functions converging to $f(z)/f(z^{p^h})$, even in seemingly simple cases such as the hypergeometric series $f_{(-1/2,1/2)}$. More deeply, one may ask whether these congruences can be used to construct the unit root of some $p$-adic cohomologies.
\end{rema}

We now state Kolchin's Theorem.

\begin{theo}[Kolchin's Theorem~\cite{kolchin}]\label{theo_Kolchin}
Let $(R,\partial)$ be a differential ring, and let $(B,\partial)$ be a differential field with $B \subset R$. Suppose that $f_1(z),\ldots, f_m(z) \in R$ satisfy $\frac{\partial f_1(z)}{f_1(z)},\ldots,\frac{\partial f_m(z)}{f_m(z)} \in B$. Then $f_1(z),\ldots, f_m(z)$ are algebraically dependent over $B$ if and only if there exist integers $a_1,\ldots, a_m \in \mathbb{Z}$, not all zero, such that
\[
f_1(z)^{a_1}\cdots f_m(z)^{a_m} \in B.
\]
\end{theo}

In the proof of Theorem~\ref{theo_alg_ind}, we apply Kolchin's Theorem with $R = \mathcal{A}_K$ and $B = E_K$.

\begin{proof}[Proof of Theorem~\ref{theo_alg_ind}]
Let $f_1(z),\ldots, f_m(z) \in \mathcal{M}\mathcal{F}(K)$. According to part~(i) of Theorem~\ref{theo_analytic_element}, $f_1(z),\ldots, f_m(z)$ belong to $1 + z\mathcal{O}_K[[z]]$, and thus $f_1(z),\ldots, f_m(z) \in \mathcal{A}_K$. It also follows from part~(iii) of Theorem~\ref{theo_analytic_element} that $f'_1(z)/f_1(z), \ldots, f'_m(z)/f_m(z)$ belong to $E_{0,K}$. Since $E_{0,K} \subset E_K \subset \mathcal{A}_K$, by Kolchin's Theorem we conclude that $f_1(z),\ldots, f_m(z)$ are algebraically dependent over $E_K$ if and only if there exist integers $a_1,\ldots, a_m \in \mathbb{Z}$, not all zero, such that
\[
f_1(z)^{a_1}\cdots f_m(z)^{a_m} \in E_K.
\]
To complete the proof, it remains to show that $h(z) = f_1(z)^{a_1} \cdots f_m(z)^{a_m} \in E_{0,K}$. For any $1 \leq i \leq m$, $f_i(z)$ converges in $D_0 = \{x \in \mathbb{C}_p : |x| < 1\}$ because $f_i(z) \in \mathcal{O}_K[[z]]$. Moreover, we have $|f_i(x)| = 1$ for all $x \in D_0$ because $f_i(0) = 1$. In particular, $|h(x)| = 1$ for all $x \in D_0$, and thus $h(z) \in E_{0,K}$.
\end{proof}

Before proving Theorem~\ref{theo_derivatives}, we establish the following lemma.

\begin{lemm}\label{lemm_derivatives}
Let $K$ be a Frobenius field, and let $f(z) \in \mathcal{M}\mathcal{F}(K)$ be such that $f^{(r)}(z) \neq 0$ for all integers $r \geq 0$. Then, for all $r \in \mathbb{N}$, $f(z)/f^{(r)}(z)$ belongs to $E_K \setminus \{0\}$.
\end{lemm}

\begin{proof}
We proceed by induction on $r \geq 1$. Since $K$ is a Frobenius field, Theorem~\ref{theo_analytic_element} implies that $\frac{f(z)}{f'(z)} \in E_K \setminus \{0\}$. Now, assume as the induction hypothesis that $\frac{f(z)}{f^{(r)}(z)} \in E_K \setminus \{0\}$. Then, the derivative $\frac{d}{dz} \left( \frac{f(z)}{f^{(r)}(z)} \right)$ also belongs to $E_K$. Clearly,
\[
\frac{d}{dz} \left( \frac{f(z)}{f^{(r)}(z)} \right) = \frac{f'(z)}{f^{(r)}(z)} - \frac{f(z) f^{(r+1)}(z)}{(f^{(r)}(z))^2}.
\]
Since, by the induction hypothesis, $f^{(r)}(z)/f(z) \in E_K$, it follows that
\[
\frac{f^{(r)}(z)}{f(z)} \cdot \frac{d}{dz}\left(\frac{f(z)}{f^{(r)}(z)}\right) = \frac{f'(z)}{f(z)} - \frac{f^{(r+1)}(z)}{f^{(r)}(z)} \in E_K.
\]
By Theorem~\ref{theo_analytic_element}, we also have $f'(z)/f(z) \in E_K$, from which it follows that $f^{(r+1)}(z)/f^{(r)}(z) \in E_K$. Consequently, we conclude that
\[
\frac{f^{(r)}(z)}{f(z)} \cdot \frac{f^{(r+1)}(z)}{f^{(r)}(z)} = \frac{f^{(r+1)}(z)}{f(z)} \in E_K \setminus \{0\}.
\]
\end{proof}

We are now ready to prove Theorem~\ref{theo_derivatives}.

\begin{proof}[Proof of Theorem~\ref{theo_derivatives}]
$(1) \Rightarrow (2)$. Let $(r_1,\ldots, r_m) \in \mathbb{N}^m$. Since $f_1(z),\ldots, f_m(z)$ are algebraically dependent over $E_K$, Theorem~\ref{theo_alg_ind} implies that there exist $a_1,\ldots, a_m \in \mathbb{Z}$, not all zero, such that
\[
f_1(z)^{a_1}\cdots f_m(z)^{a_m} \in E_{0,K} \setminus \{0\}.
\]
Furthermore, it is clear that
\begin{align*}
f_1(z)^{a_1}\cdots f_m(z)^{a_m} &= \left(\frac{f_1(z)}{f_1^{(r_1)}(z)}\right)^{a_1}\left(f_1^{(r_1)}(z)\right)^{a_1} \cdots \left(\frac{f_m(z)}{f_m^{(r_m)}(z)}\right)^{a_m}\left(f_m^{(r_m)}(z)\right)^{a_m} \\
&= \left(\frac{f_1(z)}{f_1^{(r_1)}(z)}\right)^{a_1}\cdots\left(\frac{f_m(z)}{f_m^{(r_m)}(z)}\right)^{a_m} \left(f_1^{(r_1)}(z)\right)^{a_1}\cdots\left(f_m^{(r_m)}(z)\right)^{a_m}.
\end{align*}
By Lemma~\ref{lemm_derivatives}, for all $1 \leq i \leq m$, we have $f_i(z)/f_i^{(r_i)}(z) \in E_K \setminus \{0\}$, and thus
\[
\left(\frac{f_1(z)}{f_1^{(r_1)}(z)}\right)^{a_1}\cdots\left(\frac{f_m(z)}{f_m^{(r_m)}(z)}\right)^{a_m} \in E_K \setminus \{0\}.
\]
Hence,
\[
\left(f_1^{(r_1)}(z)\right)^{a_1}\cdots\left(f_m^{(r_m)}(z)\right)^{a_m} = f_1(z)^{a_1}\cdots f_m(z)^{a_m}\left(\frac{f_1^{(r_1)}(z)}{f_1(z)}\right)^{a_1}\cdots\left(\frac{f_m^{(r_m)}(z)}{f_m(z)}\right)^{a_m} \in E_K \setminus \{0\}.
\]
Consequently, $f_1^{(r_1)},\ldots, f_m^{(r_m)}$ are algebraically dependent over $E_K$.

$(2) \Rightarrow (3)$ This direction is clear.

$(3) \Rightarrow (1)$ Suppose that there exists $(s_1,\ldots, s_m) \in \mathbb{N}^m$ such that $f_1^{(s_1)},\ldots, f_m^{(s_m)}$ are algebraically dependent over $E_K$. Then
\[
P(f_1^{(s_1)},\ldots, f_m^{(s_m)}) = \sum_{(i_1,\ldots,i_m)\in\mathbb{N}^m} b(i_1,\ldots, i_m)\left(f_1^{(s_1)}\right)^{i_1}\cdots\left(f_m^{(s_m)}\right)^{i_m} = 0,
\]
where $P(X_1,\ldots, X_m) = \sum_{(i_1,\ldots,i_m)\in\mathbb{N}^m} b(i_1,\ldots, i_m)X_1^{i_1}\cdots X_m^{i_m} \in E_K[X_1,\ldots, X_m] \setminus \{0\}$. Therefore, $Q(f_1,\ldots, f_m) = 0$, where
\[
Q(X_1,\ldots, X_m) = \sum_{(i_1,\ldots,i_m)\in\mathbb{N}^m} b(i_1,\ldots, i_m)\left(\frac{f_1^{(s_1)}}{f_1}\right)^{i_1}\cdots\left(\frac{f_m^{(s_m)}}{f_m}\right)^{i_m} X_1^{i_1}\cdots X_m^{i_m}.
\]
Since, by Lemma~\ref{lemm_derivatives}, $f_i^{(s_i)}/f_i \in E_K \setminus \{0\}$ for all $1 \leq i \leq m$, we have $Q(X_1,\ldots, X_m) \in E_K[X_1,\ldots, X_m] \setminus \{0\}$. For this reason, $f_1(z),\ldots, f_m(z)$ are algebraically dependent over $E_K$.

\end{proof}

The remainder of the paper is devoted to the proof of Theorem~\ref{theo_analytic_element}. We begin by outlining the strategy of the proof.

\section{Strategy of the proof of Theorem~\ref{theo_analytic_element}}\label{sec_stra}

In this section, we assume that $K$ is a Frobenius field with uniformizer $\pi$, and that $\mathcal{L} \in E_{0,K}[\delta]$ is MOM at zero with a strong Frobenius structure of period $h$. We also assume that $f(z) \in 1 + z K[[z]]$ is a solution to $\mathcal{L}$. Before outlining the main steps of the proof, recall that the \emph{Cartier operator} is the $\mathbb{C}_p$-linear map $\Lambda_p : \mathbb{C}_p[[z]] \to \mathbb{C}_p[[z]]$ defined by 
\[
\Lambda_p \left( \sum_{n \geq 0} a_n z^n \right) = \sum_{n \geq 0} a_{np} z^n.
\]

\paragraph{Step 1.} For every integer $m \geq 1$, we construct a monic differential operator $\mathcal{L}_m \in E_{0,K}[\delta]$ such that:
\begin{enumerate}[label=(\roman*)]
\item $\mathcal{L}_m$ is MOM at zero;
\item $\Lambda^m_p(f(z))$ is a solution to $\mathcal{L}_m$.
\end{enumerate}
In the construction of $\mathcal{L}_m$, it is crucial that $\mathcal{L}$ is MOM at zero and admits a strong Frobenius structure.

\paragraph{Step 2.} Let $A$ be the companion matrix of $\mathcal{L}$, and let $A_m$ be the companion matrix of $\mathcal{L}_m$. Then, for every integer $m \geq 1$, the matrices $A$ and $p^m A_m(z^{p^m})$ are $E_K$-equivalent. More precisely, we show that there exists a matrix $\widetilde{H}_m = (\widetilde{h}_{i,j}^{(m)})_{1 \leq i,j \leq n}$ in $\mathrm{GL}_n(E_K) \cap M_n(E_{0,K})$ such that 
\[
\delta(\widetilde{H}_m) = A\widetilde{H}_m - p^m\widetilde{H}_m A_m(z^{p^m}).
\]
Moreover, for all $1 \leq i,j \leq n$, we have $|\widetilde{h}_{i,j}^{(m)}|_{\mathcal{G}} \leq 1$ and $\widetilde{H}_m(0) = \mathrm{diag}(1, p^m, \ldots, p^{m(n-1)})$.

Steps 1 and 2 are established in Corollary~\ref{prop_ant}.

The matrices $\widetilde{H}_m$ play a central role in the proof of assertions (i) and (ii) of Theorem~\ref{theo_analytic_element}.

Let us show how to deduce assertion~(i) of Theorem~\ref{theo_analytic_element}. Write $f(z) = \sum_{j \geq 0} f_j z^j$. Recall that our goal is to prove that $f_j \in \mathcal{O}_K$ for all $j \geq 0$. As we will see, thanks to Remark~\ref{rem_frob_ac}, for every integer $m \geq 1$,
\[
\widetilde{h}^{(m)}_{1,1} \Lambda^m_p(f)(z^{p^m}) + \widetilde{h}^{(m)}_{1,2} (\delta \Lambda^m_p(f))(z^{p^m}) + \cdots + \widetilde{h}^{(m)}_{1,n} (\delta^{n-1} \Lambda^m_p f)(z^{p^m})
\]
is a solution to $\mathcal{L}$. Moreover,
\[
\sum_{i=1}^{n} \widetilde{h}^{(m)}_{1,i} (\delta^{i-1} \Lambda^m_p(f))(z^{p^m}) \in K[[z]]
\]
because $\widetilde{h}^{(m)}_{1,1}, \ldots, \widetilde{h}^{(m)}_{1,n} \in E_{0,K} \subset K[[z]]$ and clearly $\Lambda^m_p(f) \in K[[z]]$.

Since $\mathcal{L}$ is MOM at zero, as shown in Remark~\ref{rem_sol}, any solution to $\mathcal{L}$ in $K[[z]]$ is of the form $cf(z)$ for some $c \in K$. Thus, 
\[
\sum_{i=1}^{n} \widetilde{h}^{(m)}_{1,i} (\delta^{i-1} \Lambda^m_p(f))(z^{p^m}) = r f(z)
\]
for some $r \in K$. Since $f(0) = 1$ and $\widetilde{h}^{(m)}_{1,1}(0) = 1$, it follows that $r = 1$. Hence,
\begin{equation}\label{eq_str1}
\sum_{i=1}^{n} \widetilde{h}^{(m)}_{1,i} (\delta^{i-1} \Lambda^m_p(f))(z^{p^m}) = f(z).
\end{equation}
Finally, for any $2 \leq i \leq n$, we have $(\delta^{i-1} \Lambda^m_p f)(z^{p^m}) \in z^{p^m} K[[z]]$. Therefore, from \eqref{eq_str1}, we conclude that
\begin{equation}\label{eq_str}
\widetilde{h}^{(m)}_{1,1} (\Lambda^m_p(f))(z^{p^m}) \equiv f(z) \bmod{z^{p^m} K[[z]]}.
\end{equation}
Notice that $\widetilde{h}^{(m)}_{1,1} = \sum_{j \geq 0} h_j z^j \in \mathcal{O}_K[[z]]$ because $\widetilde{h}^{(m)}_{1,1} \in K[[z]]$ and $|\widetilde{h}^{(m)}_{1,1}|_{\mathcal{G}} \leq 1$. In addition, it follows from \eqref{eq_str} that $h_j = f_j$ for every $0 \leq j \leq p^m - 1$. Consequently, $f_j \in \mathcal{O}_K$ for all $0 \leq j \leq p^m - 1$. Since $m \geq 1$ is an arbitrary natural number, we conclude that $f(z) \in \mathcal{O}_K[[z]]$.

\paragraph{Step 3.} The proof of assertion~(ii) of Theorem~\ref{theo_analytic_element} proceeds in three steps:

\begin{description}[font=\normalfont\itshape]
\item[Step 3.1] In Lemma~\ref{lemm_ecuacion_nivel_1}, we prove that $f(z) / (\Lambda_p(f))(z^p)$ belongs to $E_{0,K}$. The key ingredient in the proof of this lemma is that, for all $\nu \in \{2, \ldots, n\}$,
\begin{equation}\label{eq_congruences0}
\widetilde{h}^{(m)}_{1,\nu} \equiv \widetilde{h}^{(m)}_{1,1} g_{\nu}(z^p) \bmod{\pi^m \mathcal{O}_K[[z]]},
\end{equation}
where $g_{\nu} \in E_{0,K} \cap \mathcal{O}_K[[z]]$ and $\pi$ is a uniformizer of $\mathcal{O}_K$.

The proof of the congruences \eqref{eq_congruences0} relies on Lemma~\ref{lemm_congruences}, which is technical but crucial for obtaining congruences modulo $\pi^m$ for any integer $m \geq 1$. Combining \eqref{eq_str1} and \eqref{eq_congruences0}, we obtain
\begin{equation}
\widetilde{h}^{(m)}_{1,1} d(z^p) \equiv f(z) \bmod{\pi^m \mathcal{O}_K[[z]]}, \quad \text{where} \quad d(z^p) = \sum_{\nu=1}^n g_{\nu}(z^p) (\delta^{\nu-1} \Lambda^m_p(f))(z^{p^m}).
\end{equation}
Applying $\Lambda_p$ to this equality yields $d(z) \equiv \Lambda_p(f) / \Lambda_p(\widetilde{h}^{(m)}_{1,1}) \bmod{\pi^m \mathcal{O}_K[[z]]}$.\footnote{Here we use the fact that $\Lambda_p(h(z)t(z^p)) = \Lambda_p(h(z))t(z)$ for any $h(z), t(z) \in \mathbb{C}_p[[z]]$.} Hence, we deduce that
\[
\frac{\widetilde{h}^{(m)}_{1,1}}{(\Lambda_p(\widetilde{h}^{(m)}_{1,1}))(z^p)} \equiv \frac{f(z)}{(\Lambda_p(f))(z^p)} \bmod{\pi^m \mathcal{O}_K[[z]]}.
\]
As we will show, $\frac{\widetilde{h}^{(m)}_{1,1}}{(\Lambda_p(\widetilde{h}^{(m)}_{1,1}))(z^p)}$ belongs to $E_{0,K}$. Since $m \geq 1$ is arbitrary, we conclude that $f(z) / (\Lambda_p(f))(z^p) \in E_{0,K}$.

In a similar manner, Corollary~\ref{prop_induction} shows that, for all integers $k \geq 1$, the quotient $\frac{f(z)}{(\Lambda^k_p f)(z^{p^k})}$ belongs to $E_{0,K}$.\smallskip

\item[Step 3.2] In Proposition~\ref{prop_rational}, we prove that $\frac{f(z)}{\Lambda^{kh}_p f(z)} \in E_{0,K}$ for any integer $k \geq 1$, where $h$ is the period of the strong Frobenius structure.\smallskip

\item[Step 3.3.] We deduce that $f(z)/f(z^{p^h})$ belongs to $E_{0,K}$ as follows. Let $k \geq 1$ be an integer. By Step 3.1, $\frac{f(z)}{(\Lambda^{kh}_p f)(z^{p^{kh}})}$ belongs to $E_{0,K}$ for all $k \geq 1$. In particular, there exists $P_{kh}(z) \in K_0(z)$ such that
\[
P_{kh}(z)(\Lambda_p^{kh}(f))(z^{p^{kh}}) \equiv f(z) \bmod{\pi^{kh} \mathcal{O}_K[[z]]}.
\]
Rewriting this, we obtain
\[
P_{kh}(z) \frac{(\Lambda_p^{kh}(f))(z^{p^{kh}})}{f(z^{p^{kh}})} f(z^{p^{kh}}) \equiv f(z) \bmod{\pi^{kh} \mathcal{O}_K[[z]]}.
\]
By Step 3.2, there exists $Q_k(z) \in K_0(z)$ with norm $1$ such that
\[
Q_k(z) \equiv \frac{\Lambda_p^{kh}(f(z))}{f(z)} \bmod{\pi^{kh} \mathcal{O}_K[[z]]}.
\]
Combining these congruences yields
\[
P_{kh}(z) Q_k(z^{p^{kh}}) f(z^{p^{kh}}) \equiv f(z) \bmod{\pi^{kh} \mathcal{O}_K[[z]]}.
\]
Setting $B_k(z) = P_{kh}(z) Q_k(z^{p^{kh}})$, we obtain
\[
f(z) \equiv B_k(z) f(z^{p^{kh}}) \bmod{\pi^{kh} \mathcal{O}_K[[z]]}.
\]
Hence,
\[
\frac{f(z)}{f(z^{p^h})} \equiv \frac{B_{k+1}(z)}{B_k(z^{p^h})} \bmod{\pi^{kh} \mathcal{O}_K[[z]]}.
\]
Since $k \geq 1$ is arbitrary, this congruence implies that $f(z)/f(z^{p^h}) \in E_{0,K}$.
\end{description}

\paragraph{Step 4.} Finally, to prove assertion~(iii) of Theorem~\ref{theo_analytic_element}, we define, for every integer $m \geq 1$, 
\[
H_m(z) = H(z) \cdots H(z^{p^{m-1}}), \quad \text{where } H(z) = \frac{f(z)}{f(z^{p^h})}.
\]
Notice that $H_m \in E_{0,K}$ because $H \in E_{0,K}$ by Step 3. We then prove that
\[
\frac{f'(z)}{f(z)} = \lim_{m \to \infty} \frac{H'_m(z)}{H_m(z)}.
\]
Consequently, $f'(z)/f(z)$ is an element of $E_{0,K}$.

\section{Constructing Frobenius antecedents and Cartier operators}\label{sec_frob_ant}

In this section, we prove Steps~1 and 2 described in Section~\ref{sec_stra}, which, as explained there, lead to the proof of assertion~(i) of Theorem~\ref{theo_analytic_element}. To this end, we provide an explicit construction of the \emph{Frobenius antecedent} for MOM differential operators $\mathcal{L} \in E_p[\delta]$. Our construction is based on the work of Christol~\cite{christolpadique}. 

Recall that a matrix $A \in M_n(E_p)$ is said to admit a Frobenius antecedent if there exists a matrix $B \in M_n(E_p)$ such that $A$ and $pB(z^p)$ are $E_p$-equivalent. The main goal of this section is to prove that if $\mathcal{L} \in E_p[\delta]$ is a MOM differential operator equipped with a strong Frobenius structure, and $A$ is its companion matrix, then there exists another MOM differential operator $\mathcal{L}_1 \in E_p[\delta]$, also equipped with a strong Frobenius structure, such that $A$ and $pA_1(z^p)$ are $E_p$-equivalent, where $A_1$ is the companion matrix of $\mathcal{L}_1$. By definition, this means that there exists a matrix $\widetilde{H}_1 \in \mathrm{GL}_n(E_p)$ such that 
\[
\delta(\widetilde{H}_1) = A\widetilde{H}_1 - p\widetilde{H}_1 A_1(z^p).
\]
In Proposition~\ref{lemm_paso_base}, using the \emph{Cartier operator}, we provide an explicit construction of $\widetilde{H}_1$ and establish certain $p$-adic properties of $\widetilde{H}_1$ that play a crucial role in the proof of Theorem~\ref{theo_analytic_element}. For this purpose, we first introduce the set $\mathcal{M}_K$.

\begin{defi}\label{def_m_p}
Let $K$ be a complete extension of $\mathbb{Q}_p$. The set $\mathcal{M}_K$ consists of square matrices $A$ with entries in $E_K$ satisfying the following conditions:
\begin{enumerate}[label=(\roman*)]
\item $A$ has entries in $E_{0,K}$.
\item The eigenvalues of $A(0)$ are all equal to zero.
\item There exists an invertible matrix $U$ with entries in the generic disc $D(t,1^-)$ such that $\delta U = AU$.
\end{enumerate}
\end{defi}

The relation between the sets $\mathcal{M}\mathcal{F}(K)$ and $\mathcal{M}_K$ is given by the following observation: if $\mathcal{L} \in E_{0,K}[\delta]$ is MOM at zero and admits a strong Frobenius structure, then its companion matrix $A$ belongs to $\mathcal{M}_K$. Indeed, if $A$ is the companion matrix of $\mathcal{L}$, it is clear that $A$ has entries in $E_{0,K}$, and the MOM condition is equivalent to requiring that all eigenvalues of $A(0)$ are zero. Finally, since $\mathcal{L}$ has a strong Frobenius structure, it follows from Propositions~4.1.2, 4.6.4, and 4.7.2 of \cite{C83} that there exists an invertible matrix $U$ with entries in the generic disc $D(t,1^-)$ such that $\delta U = AU$.

\begin{rema}\label{rem_sol}(\textbf{Fundamental matrix of Solutions})

Let $K$ be a complete extension of $\mathbb{Q}_p$. Suppose that $A \in M_n(E_{0,K})$ and that all eigenvalues of $A(0)$ are zero. In particular, the eigenvalues of $A(0)$ are prepared.\footnote{Given a matrix $C \in M_n(\mathbb{C}_p)$, we say that its eigenvalues are \emph{prepared} if the following condition holds: for any two eigenvalues $\lambda$ and $\alpha$, $\lambda - \alpha$ is an integer if and only if $\lambda = \alpha$.} Since $E_{0,K} \subset K[[z]]$, we have $A \in M_n(K[[z]])$. It then follows from~\cite[Chap.~III, Proposition~8.5]{Dworkgfunciones} that the differential system $\delta X = AX$ admits a \textit{fundamental matrix of solutions} of the form $Y_A z^{A(0)}$, where $Y_A \in \mathrm{GL}_n(K[[z]])$ and $Y_A(0) = I_n$. Recall that
\[
z^{A(0)} = \sum_{j \ge 0} A(0)^j \frac{(\log z)^j}{j!}.
\]
Since all eigenvalues of $A(0)$ are zero, $A(0)$ is nilpotent with $A(0)^n = 0$, and therefore
\[
z^{A(0)} = \sum_{j=0}^{n-1} A(0)^j \frac{(\log z)^j}{j!}.
\]
The matrix $Y_A$ is referred to as the \emph{uniform part} of the system $\delta X = AX$.

\smallskip

Now, let $L \in E_{0,K}[\delta]$ and let $A$ be its companion matrix. If $L$ is MOM at zero, then by definition all eigenvalues of $A(0)$ are zero. Hence, by the preceding argument, $L$ admits a solution $f(z) \in 1 + z K[[z]]$. Moreover, since $\log z$ is transcendental over $K[[z]]$, any solution $g(z) \in K[[z]]$ to $L$ is of the form $g(z) = c f(z)$ for some constant $c \in K$.
\end{rema}
 
 \subsection{Constructing Frobenius antecedents}
 
An important ingredient in the construction of the Frobenius antecedent for elements of $\mathcal{M}_K$ is the \emph{Cartier operator}. Recall that the Cartier operator is the $\mathbb{C}_p$-linear map $\Lambda_p : \mathbb{C}_p[[z]] \to \mathbb{C}_p[[z]]$ defined by 
\[
\Lambda_p \left( \sum_{n \ge 0} a_n z^n \right) = \sum_{n \ge 0} a_{np} z^n.
\]
If $A$ is a matrix with entries in $\mathbb{C}_p[[z]]$, then $\Lambda_p(A)$ denotes the matrix obtained by applying $\Lambda_p$ entrywise to $A$.

In~\cite{christolpadique}, Christol provided a construction of the Frobenius antecedent for a set more general than $\mathcal{M}_K$. However, our key observation is that when $A$ belongs to $\mathcal{M}_K$, this construction can be made more explicit. This is carried out in Proposition~\ref{lemm_paso_base}. To prove this proposition, we need the following lemma, whose proof follows the same lines as that of~\cite[Lemma~5.1]{christolpadique}.

Finally, given a matrix $A = (a_{i,j})_{1 \le i,j \le n}$ with entries in $E_K$, we define its norm by
\[
\|A\| := \max_{1 \le i,j \le n} \{ |a_{i,j}|_{\mathcal{G}} \}.
\]
  
  \begin{lemm}\label{lemm_paso_0}
Let $K$ be a finite extension of $\mathbb{Q}_p$. Let $\mathcal{L} \in E_{0,K}[\delta]$ be a MOM differential operator of order $n$, let $A$ be the companion matrix of $\mathcal{L}$, and let $Y_A z^{A(0)}$ be a fundamental matrix of solutions to $\delta X = AX$, where $Y_A$ is the uniform part. If $A$ belongs to $\mathcal{M}_K$ and $\|A\| \le 1$, then there exists a matrix $F \in M_n(E_{0,K})$ such that $A$ and $p F(z^p)$ are $E_K$-equivalent. Moreover, the matrix 
\[
H_0 := Y_A (\Lambda_p(Y_A)(z^p))^{-1}
\]
belongs to $\mathrm{GL}_n(E_{0,K})$, satisfies $\|H_0\| = 1 = \|H_0^{-1}\|$, and fulfills
\[
\delta H_0 = p F(z^p) H_0 - H_0 A.
\]
In addition, the matrix $\Lambda_p(Y_A) z^{\frac{1}{p} A(0)}$ is a fundamental matrix of solutions to $\delta X = FX$.
\end{lemm}

\begin{proof}
We first construct the matrix $H_0$. To this end, consider the sequence of matrices $\{A_j(z)\}_{j \ge 0}$ defined by $A_0(z) = I_n$ and
\[
A_{j+1}(z) = \delta A_j(z) + A_j(z)(A(z) - j I_n) \quad \text{for } j \ge 0,
\]
where $I_n$ denotes the identity matrix. By assumption, $A \in M_n(E_{0,K})$, and it follows inductively that $A_j(z) \in M_n(E_{0,K})$ for all $j \ge 0$.

Next, recall that $X^p - 1 = (X - 1)\bigl((X - 1)^{p-1} + p\, t(X)\bigr)$ with $t(X) \in \mathbb{Z}[X]$. Therefore, if $\xi^p = 1$, then $|\xi - 1| \le |p|^{1/(p-1)}$. Moreover, for all integers $j \ge 1$, we have $|p|^{j/(p-1)} < |j!|$, and hence
\[
\left| \frac{(\xi - 1)^j}{j!} \right| < 1.
\]
Since the $p$-adic norm is non-Archimedean and $\sum_{\xi^p = 1} (\xi - 1)^j / j!$ is a rational number, we conclude that for all $j \ge 0$,
\[
\left| \sum_{\xi^p = 1} \frac{(\xi - 1)^j}{j!} \right| \le \frac{1}{p}.
\]
Consequently, $\left| \sum_{\xi^p = 1} \frac{(\xi - 1)^j}{p j!} \right| \le 1$ for all integers $j \ge 0$.

We then define
\[
H_0 = \sum_{j \ge 0} A_j(z) \left( \sum_{\xi^p = 1} \frac{(\xi - 1)^j}{p j!} \right).
\]

We now show that $H_0 \in \mathrm{GL}_n(E_{0,K})$. First, we prove that $H_0 \in M_n(E_{0,K})$, which is equivalent to showing that
\begin{equation}\label{eq_lim}
\lim_{j \to \infty} \|A_j\| \left| \sum_{\xi^p = 1} \frac{(\xi - 1)^j}{p j!} \right| = 0,
\end{equation}
since $E_{0,K}$ is complete with respect to the non-Archimedean Gauss norm. Indeed, because $A \in \mathcal{M}_K$, there exists an invertible matrix $U$ on the generic disc $D(t,1^-)$ such that $\delta U = A U$. Thus, it follows from Remark~b) of \cite[p.~155]{christolpadique} that $\lim_{j \to \infty} \|A_j\| = 0$. Therefore, \eqref{eq_lim} follows immediately since $\left| \sum_{\xi^p = 1} \frac{(\xi - 1)^j}{p j!} \right| \le 1$ for all $j \ge 0$.

Furthermore, we have $\|H_0\| \le 1$. Indeed, by assumption $\|A\| \le 1$, and consequently $\|A_j\| \le 1$ for all $j \ge 0$. We have also shown that $\left| \sum_{\xi^p = 1} \frac{(\xi - 1)^j}{p j!} \right| \le 1$ for all $j \ge 0$, and therefore $\|A_j\| \left| \sum_{\xi^p = 1} \frac{(\xi - 1)^j}{p j!} \right| \le 1$. As the Gauss norm is non-Archimedean, we conclude that $\|H_0\| \le 1$.

We now proceed to prove that $H_0 = \Lambda_p(Y_A)(z^p) Y_A^{-1}$. Write $Y_A = \sum_{j \ge 0} Y_j z^j$. Since $A \in \mathcal{M}_K$, we have $A(0)^n = 0$. It then follows from \cite[p.~165]{christolpadique} that $H_0 Y_A = \sum_{j \ge 0} Y_{jp} z^{jp} = \Lambda_p(Y_A)(z^p)$. Thus, $H_0 = \Lambda_p(Y_A)(z^p) Y_A^{-1}$. Finally, $H_0 \in \mathrm{GL}_n(E_{0,K})$ because $H_0(0) = I_n$. In particular, $\|H_0\| = 1 = \|H_0^{-1}\|$.

Next, we construct the matrix $F$. To this end, we set
\begin{equation}\label{eq_frob_m'}
F(z) = \left[ \delta(\Lambda_p(Y_A)) + \frac{1}{p} \Lambda_p(Y_A) A(0) \right] (\Lambda_p(Y_A))^{-1}.
\end{equation}
Since $H_0 Y_A = \Lambda_p(Y_A)(z^p)$, Equation~\eqref{eq_frob_m'} yields
\[
p F(z^p) = \bigl[ p(\delta(\Lambda_p(Y_A)))(z^p) + H_0 Y_A A(0) \bigr] \bigl[ Y_A^{-1} H_0^{-1} \bigr].
\]
We also have
\[
\delta(H_0) Y_A + H_0(\delta Y_A) = \delta(H_0 Y_A) = \delta(\Lambda_p(Y_A)(z^p)) = p(\delta(\Lambda_p(Y_A)))(z^p).
\]
Since $\delta(Y_A z^{A(0)}) = A Y_A z^{A(0)}$, it follows that $\delta Y_A = A Y_A - Y_A A(0)$. Thus,
\[
p(\delta(\Lambda_p(Y_A)))(z^p) = \delta(H_0) Y_A + H_0 \bigl[ A Y_A - Y_A A(0) \bigr].
\]
Then,
\begin{align*}
p F(z^p) &= \bigl[ p(\delta(\Lambda_p(Y_A)))(z^p) + H_0 Y_A A(0) \bigr] [Y_A^{-1} H_0^{-1}] \\
&= \bigl[ (\delta H_0) Y_A + H_0 A Y_A \bigr] [Y_A^{-1} H_0^{-1}] \\
&= (\delta H_0) H_0^{-1} + H_0 A H_0^{-1}.
\end{align*}
Consequently, $\delta H_0 = p F(z^p) H_0 - H_0 A$. Since $p F(z^p) = (\delta H_0 + H_0 A) H_0^{-1}$ with $H_0 \in \mathrm{GL}_n(E_{0,K})$ and $A \in M_n(E_{0,K})$, we obtain $F(z) \in M_n(E_{0,K})$.

Finally, we show that the matrix $\Lambda_p(Y_A) z^{\frac{1}{p} A(0)}$ is a fundamental matrix of solutions to $\delta X = FX$. Clearly, $\Lambda_p(Y_A)(0) = I_n$, and from Equality~\eqref{eq_frob_m'} we obtain $F(z) \Lambda_p(Y_A) = \delta(\Lambda_p(Y_A)) + \Lambda_p(Y_A) \frac{1}{p} A(0)$. Thus,
\begin{align*}
\delta \!\left( \Lambda_p(Y_A) z^{\frac{1}{p} A(0)} \right) &= \delta(\Lambda_p(Y_A)) z^{\frac{1}{p} A(0)} + \Lambda_p(Y_A) \frac{1}{p} A(0) z^{\frac{1}{p} A(0)} \\
&= \left( \delta(\Lambda_p(Y_A)) + \Lambda_p(Y_A) \frac{1}{p} A(0) \right) z^{\frac{1}{p} A(0)} \\
&= F(z) \Lambda_p(Y_A) z^{\frac{1}{p} A(0)}.
\end{align*}
This completes the proof.
\end{proof}

The key new ingredient in our construction of the Frobenius antecedent is provided by the following proposition. Indeed, it follows from this proposition that if $\mathcal{L} \in E_{0,K}[\delta]$ is a MOM differential operator equipped with a strong Frobenius structure and $A$ is its companion matrix, then there exists a MOM differential operator $\mathcal{L}_1 \in E_{0,K}[\delta]$, also equipped with a strong Frobenius structure, such that $A$ and $p A_1(z^p)$ are $E_p$-equivalent, where $A_1$ denotes the companion matrix of $\mathcal{L}_1$. The construction of $\mathcal{L}_1$ relies on the matrix $F$ introduced in Lemma~\ref{lemm_paso_0}.

\begin{prop}\label{lemm_paso_base}
Let the assumptions be as in Lemma~\ref{lemm_paso_0}, and let $f(z) \in 1 + z K[[z]]$ be a solution to $\mathcal{L}$. Then there exists a MOM differential operator $\mathcal{L}_1 \in E_{0,K}[\delta]$ of order $n$ such that:
\begin{enumerate}[label=(\alph*)]
  \item $\Lambda_p(f)$ is a solution to $\mathcal{L}_1$.
  
  \item If $A_1$ denotes the companion matrix of $\mathcal{L}_1$, then a fundamental matrix of solutions to $\delta X = A_1 X$ is given by
  \[
  \mathrm{diag}\left(1, \frac{1}{p}, \ldots, \frac{1}{p^{n-1}}\right) \Lambda_p(Y_A) \, \mathrm{diag}(1, p, \ldots, p^{n-1}) \, z^{A(0)}.
  \]
  
  \item $A_1$ belongs to $\mathcal{M}_K$ and satisfies $\|A_1\| \le 1$.
  
  \item The matrices $A$ and $p A_1(z^p)$ are $E_K$-equivalent. Moreover,
  \[
  \delta \widetilde{H}_1 = A \widetilde{H}_1 - p \widetilde{H}_1 A_1(z^p),
  \]
  where
  \[
  \widetilde{H}_1 := Y_A (\Lambda_p(Y_A)(z^p))^{-1} \mathrm{diag}(1, p, \ldots, p^{n-1}) \in \mathrm{GL}_n(E_{0,K}),
  \]
  and $\|\widetilde{H}_1\| = 1$.
\end{enumerate}
\end{prop}
 
\begin{proof}
We begin by constructing the MOM differential operator $\mathcal{L}_1 \in E_{0,K}[\delta]$ of order $n$. By Lemma~\ref{lemm_paso_0}, there exists a matrix $F \in M_n(E_{0,K})$ such that $A$ and $p F(z^p)$ are $E_K$-equivalent. Write $F = (b_{i,j})_{1 \le i,j \le n}$. We then define
\[
\mathcal{L}_1 := \delta^n - b_{n,n}\, \delta^{n-1} - \frac{1}{p} b_{n,n-1} \, \delta^{n-2} - \cdots - \frac{1}{p^{n-i}} b_{n,i} \, \delta^{i-1} - \cdots - \frac{1}{p^{n-2}} b_{n,2} \, \delta - \frac{1}{p^{n-1}} b_{n,1}.
\]

Clearly, $\mathcal{L}_1 \in E_{0,K}[\delta]$ because the entries of $F$ lie in $E_{0,K}$, and its order is evidently $n$. We now prove that $\mathcal{L}_1$ is MOM at zero. Indeed, from Lemma~\ref{lemm_paso_0}, we have
\begin{equation}\label{eq_h_0}
\delta H_0 = p F(z^p) H_0 - H_0 A,
\end{equation}
where $H_0 = Y_A (\Lambda_p(Y_A)(z^p))^{-1} \in \mathrm{GL}_n(E_{0,K})$. Since $Y_A(0) = I_n$ (as $Y_A$ is the uniform part of $\delta X = A X$), it follows that $H_0(0) = I_n$. Substituting $z = 0$ into \eqref{eq_h_0} gives $F(0) = \frac{1}{p} A(0)$. Consequently, for each $i \in \{1, \ldots, n\}$, we have $b_{n,i}(0) = \frac{1}{p} a_{n,i}(0)$. By assumption, $A \in \mathcal{M}_K$, so all eigenvalues of $A(0)$ are zero. Since $A$ is the companion matrix of $\mathcal{L}$, this implies $a_{n,i}(0) = 0$ for all $1 \le i \le n$. Therefore, $b_{n,i}(0) = 0$ for all $1 \le i \le n$.

We now proceed to verify statements (a)–(d).

\smallskip
\noindent
\emph{Proof of (a)}. We show that $\Lambda_p(f)$ is a solution to $\mathcal{L}_1$. By Lemma~\ref{lemm_paso_0}, a fundamental matrix of solutions to $\delta X = F X$ is given by $\Lambda_p(Y_A) z^{\frac{1}{p} A(0)}$. Write $Y_A = (f_{i,j})_{1 \le i,j \le n}$. Since $Y_A$ is the uniform part of $\delta X = A X$ and $f(z) \in 1 + z K[[z]]$ is a solution to $\mathcal{L}$, we have $f_{i,1} = \delta^{i-1} f$ for all $i \in \{1,\dots,n\}$. Hence, the vector 
\[
\begin{pmatrix} \Lambda_p(f_{1,1}) \\ \Lambda_p(f_{2,1}) \\ \vdots \\ \Lambda_p(f_{n,1}) \end{pmatrix}
=
\begin{pmatrix} \Lambda_p(f) \\ \Lambda_p(\delta f) \\ \vdots \\ \Lambda_p(\delta^{n-1} f) \end{pmatrix}
\]
is a solution to $\delta \vec{y} = F \vec{y}$. Therefore,
\[
F
\begin{pmatrix} 
\Lambda_p(f) \\ \Lambda_p(\delta f) \\ \vdots \\ \Lambda_p(\delta^{n-1} f) 
\end{pmatrix}
=
\begin{pmatrix} 
\delta(\Lambda_p(f)) \\ \delta(\Lambda_p(\delta f)) \\ \vdots \\ \delta(\Lambda_p(\delta^{n-1} f)) 
\end{pmatrix}.
\]
Expanding the last row of this system, we obtain
\[
b_{n,1} \Lambda_p(f) + b_{n,2} \Lambda_p(\delta f) + \cdots + b_{n,n} \Lambda_p(\delta^{n-1} f) = \delta(\Lambda_p(\delta^{n-1} f)).
\]
Using the identity $\Lambda_p \circ \delta = p \, \delta \circ \Lambda_p$, we then have
\[
b_{n,1} \Lambda_p(f) + p b_{n,2} \delta(\Lambda_p(f)) + \cdots + p^{i-1} b_{n,i} \delta^{i-1}(\Lambda_p(f)) + \cdots + p^{n-1} b_{n,n} \delta^{n-1}(\Lambda_p(f)) = p^{n-1} \delta^n(\Lambda_p(f)).
\]
Dividing both sides by $p^{n-1}$, we conclude that $\Lambda_p(f)$ is a solution to $\mathcal{L}_1$, as claimed.

\smallskip
\noindent
\emph{Proof of (b)}. Recall that $Y_A = (f_{i,j})_{1 \le i,j \le n}$. We first show that, for all integers $i \in \{1, \ldots, n-1\}$ and $k \in \{1, \ldots, n\}$, 
\[
f_{i,k-1} + \delta f_{i,k} = f_{i+1,k}.
\]
To this end, define
\[
\phi_{i,j} := \sum_{k=1}^j f_{i,k} \frac{(\log z)^{j-k}}{(j-k)!} \quad \text{for all } i,j \in \{1, \ldots, n\}.
\]
Notice that $Y_A z^{A(0)} = (\phi_{i,j})_{1 \le i,j \le n}$. Since $Y_A z^{A(0)}$ is a fundamental matrix of $\delta X = A X$ and $A$ is the companion matrix of $\mathcal{L}$, it follows that for every $i \in \{2, \ldots, n\}$ and $j \in \{1, \ldots, n\}$, we have $\phi_{i,j} = \delta \phi_{i-1,j}$. Therefore,
\[
\sum_{k=1}^j f_{i+1,k} \frac{(\log z)^{j-k}}{(j-k)!} = \phi_{i+1,j} = \delta\phi_{i,j} = \delta(f_{i,1}) \frac{(\log z)^{j-1}}{(j-1)!} + \sum_{k=2}^j \bigl(f_{i,k-1} + \delta(f_{i,k})\bigr) \frac{(\log z)^{j-k}}{(j-k)!}.
\]
Since $\log z$ is transcendental over $K[[z]]$, comparing coefficients of $(\log z)^{j-k}$ in the above equality yields
\[
f_{i,k-1} + \delta f_{i,k} = f_{i+1,k} \quad \text{for all } k \in \{1, \ldots, n\}.
\]
We now prove that 
\[
T = \mathrm{diag}\left(1, \frac{1}{p}, \ldots, \frac{1}{p^{n-1}}\right) \Lambda_p(Y_A) \, \mathrm{diag}(1, p, \ldots, p^{n-1}) \, z^{A(0)}
\]
is a fundamental matrix of solutions to $\delta X = A_1 X$, where $A_1$ is the companion matrix of $\mathcal{L}_1$.

To this end, we first show that, for every $j \in \{1,\ldots,n\}$, 
\[
\theta_j := \sum_{k=1}^{j} p^{k-1} \Lambda_p(f_{1,k}) \frac{(\log z)^{j-k}}{(j-k)!}
\]
is a solution to $\mathcal{L}_1$. Notice that the vector $(\theta_1, \theta_2, \ldots, \theta_n)$ corresponds to the first row of $T$.

Next, set $\Lambda_p(Y_A) z^{\frac{1}{p} A(0)} = (\eta_{i,j})_{1 \le i,j \le n}$. Then, for every $i,j \in \{1, \ldots, n\}$,
\[
\eta_{i,j} = \sum_{k=1}^{j} \Lambda_p(f_{i,k}) \frac{(\log z)^{j-k}}{p^{j-k} (j-k)!}.
\]
We claim that, for every $j,l \in \{1,\ldots,n\}$,
\[
\frac{1}{p^{n-l}} \delta^{l-1} \theta_j = \frac{1}{p^{n-j}} \eta_{l,j}.
\]
We prove this by induction on $l$. For $l=1$, it is immediate from the definitions of $\theta_j$ and $\eta_{1,j}$ that
\[
\frac{1}{p^{n-1}} \theta_j = \frac{1}{p^{n-j}} \eta_{1,j}.
\]
Assume now that for some $l \in \{1,\ldots,n-1\}$, we have $\frac{1}{p^{n-l}} \delta^{l-1} \theta_j = \frac{1}{p^{n-j}} \eta_{l,j}$. Differentiating both sides yields $\frac{1}{p^{n-l}} \delta^l \theta_j = \frac{1}{p^{n-j}} \delta(\eta_{l,j})$. Since $\Lambda_p \circ \delta = p \delta \circ \Lambda_p$ and $f_{l,k-1} + \delta(f_{l,k}) = f_{l+1,k}$ for all $k \in \{1,\ldots, n\}$, we have
\begin{align*}
\delta(\eta_{l,j}) &= \sum_{k=1}^j \delta(\Lambda_p(f_{l,k})) \frac{(\log z)^{j-k}}{p^{j-k}(j-k)!} + \Lambda_p(f_{l,k}) \frac{(\log z)^{j-k-1}}{p^{j-k}(j-k-1)!} \\
&= \sum_{k=1}^j \frac{1}{p} \Lambda_p(\delta(f_{l,k})) \frac{(\log z)^{j-k}}{p^{j-k}(j-k)!} + \frac{1}{p} \Lambda_p(f_{l,k}) \frac{(\log z)^{j-k-1}}{p^{j-k-1}(j-k-1)!} \\
&= \frac{1}{p} \left[ \Lambda_p(\delta(f_{l,1})) \frac{(\log z)^{j-1}}{p^{j-1}(j-1)!} + \sum_{k=2}^j \Lambda_p(f_{l,k-1} + \delta f_{l,k}) \frac{(\log z)^{j-k}}{p^{j-k}(j-k)!} \right] \\
&= \frac{1}{p} \left[ \Lambda_p(f_{l+1,1}) \frac{(\log z)^{j-1}}{p^{j-1}(j-1)!} + \sum_{k=2}^j \Lambda_p(f_{l+1,k}) \frac{(\log z)^{j-k}}{p^{j-k}(j-k)!} \right] \\
&= \frac{1}{p} \eta_{l+1,j}.
\end{align*}
Thus, $\frac{1}{p^{n-l}} \delta^l \theta_j = \frac{1}{p^{n-j}} \delta(\eta_{l,j}) = \frac{1}{p^{n-j}} \left( \frac{1}{p} \eta_{l+1,j} \right)$, whence
\[
\frac{1}{p^{n-l-1}} \delta^l \theta_j = \frac{1}{p^{n-j}} \eta_{l+1,j}.
\]
This completes the induction, establishing that $\frac{1}{p^{n-l}} \delta^{l-1} \theta_j = \frac{1}{p^{n-j}} \eta_{l,j}$ for all $j,l \in \{1,\ldots,n\}$.

Since $(\eta_{i,j})_{1 \le i,j \le n} = \Lambda_p(Y_A) z^{\frac{1}{p} A(0)}$ is a fundamental matrix of $\delta X = F X$, it follows that for every $j \in \{1,\ldots,n\}$,
\[
F
\begin{pmatrix} 
\eta_{1,j} \\
\eta_{2,j} \\
\vdots \\
\eta_{n,j}
\end{pmatrix}
=
\begin{pmatrix} 
\delta(\eta_{1,j}) \\
\delta(\eta_{2,j}) \\
\vdots \\
\delta(\eta_{n,j})
\end{pmatrix}.
\]
Thus, for every $j \in \{1,\ldots,n\}$,
\[
b_{n,1} \eta_{1,j} + b_{n,2} \eta_{2,j} + \cdots + b_{n,k} \eta_{k,j} + \cdots + b_{n,n} \eta_{n,j} = \delta(\eta_{n,j}).
\]
Dividing by $p^{n-j}$, we get
\[
b_{n,1} \frac{1}{p^{n-j}} \eta_{1,j} + b_{n,2} \frac{1}{p^{n-j}} \eta_{2,j} + \cdots + b_{n,k} \frac{1}{p^{n-j}} \eta_{k,j} + \cdots + b_{n,n} \frac{1}{p^{n-j}} \eta_{n,j} = \frac{1}{p^{n-j}} \delta(\eta_{n,j}).
\]
Using $\frac{1}{p^{n-l}} \delta^{l-1} \theta_j = \frac{1}{p^{n-j}} \eta_{l,j}$, this becomes
\[
\frac{1}{p^{n-1}} b_{n,1} \theta_j + \frac{1}{p^{n-2}} b_{n,2} \delta \theta_j + \cdots + \frac{1}{p^{n-k}} b_{n,k} \delta^{k-1} \theta_j + \cdots + b_{n,n} \delta^{n-1} \theta_j = \delta^n \theta_j.
\]
Therefore, $\theta_j$ is a solution to $\mathcal{L}_1$.

Moreover, $\theta_1, \ldots, \theta_n$ are linearly independent over $K$ because $\log z$ is transcendental over $K[[z]]$. Since $A_1$ is the companion matrix of $\mathcal{L}_1$, it follows that the matrix $(\delta^{i-1} \theta_j)_{1 \le i,j \le n}$ is a fundamental matrix of solutions to $\delta X = A_1 X$. We now show that $T = (\delta^{i-1} \theta_j)_{1 \le i,j \le n}$. For every $i,j \in \{1, \ldots, n\}$, we have $\delta^{i-1} \theta_j = \frac{p^j}{p^i} \eta_{i,j}$. Hence,
\[
\delta^{i-1} \theta_j = \frac{p^j}{p^i} \sum_{k=1}^j \Lambda_p(f_{i,k}) \frac{(\log z)^{j-k}}{p^{j-k}(j-k)!} = \sum_{k=1}^j \frac{\Lambda_p(f_{i,k})}{p^{i-k}} \frac{(\log z)^{j-k}}{(j-k)!}.
\]
Since
\[
\mathrm{diag}\left(1, \frac{1}{p}, \ldots, \frac{1}{p^{n-1}}\right) \Lambda_p(Y_A) \, \mathrm{diag}(1, p, \ldots, p^{n-1}) = \left( \frac{\Lambda_p(f_{i,k})}{p^{i-k}} \right)_{1 \le i,k \le n},
\]
we conclude that $T = (\delta^{i-1} \theta_j)_{1 \le i,j \le n}$, and thus $T$ is indeed a fundamental matrix of solutions to $\delta X = A_1 X$.

\smallskip
\noindent
\emph{Proof of (c)}. Since $\mathcal{L}_1 \in E_{0,K}[\delta]$, $A_1$ is a matrix with entries in $E_{0,K}$. Furthermore, we have already established that $b_{n,i}(0) = 0$ for all $1 \le i \le n$. Thus, all eigenvalues of $A_1(0)$ are equal to zero.

We now prove that there exists an invertible matrix $U_1$ with entries in $D(t,1^-)$ such that $\delta U_1 = A_1 U_1$. Recall that $Y_A = (f_{i,j})_{1 \le i,j \le n}$ is the uniform part of $\delta X = A X$. Since there exists an invertible matrix $U$ with entries in $D(t,1^-)$ such that $\delta U = A U$, it follows from Theorem~2 of \cite{christolpadique} that for all $i,j \in \{1, \ldots, n\}$, the radius of convergence of $f_{i,j}$ is at least $1$. Because the $p$-adic norm is non-Archimedean, the radius of convergence of $\Lambda_p(f_{i,j})$ is also at least $1$ for all $i,j \in \{1, \ldots, n\}$. In particular, each entry of
\[
Y_{A_1} = \mathrm{diag}\left(1, \frac{1}{p}, \ldots, \frac{1}{p^{n-1}}\right) \Lambda_p(Y_A) \, \mathrm{diag}(1, p, \ldots, p^{n-1})
\]
has a radius of convergence at least $1$. Moreover, statement (b) implies that $Y_{A_1}$ is the uniform part of $\delta X = A_1 X$. Applying Theorem~2 of \cite{christolpadique} again, we conclude that there exists an invertible matrix $U_1$ with entries in $D(t,1^-)$ such that $\delta U_1 = A_1 U_1$. Consequently, $A_1 \in \mathcal{M}_K$ and,  by the Frobenius--Dwork Theorem \cite[Proposition~8.1]{gillesff}, we  also conclude that $\|A_1\| \le 1$.

\smallskip
\noindent
\emph{Proof of (d)}.We define
\[
\widetilde{H}_1 := Y_A (\Lambda_p(Y_A)(z^p))^{-1} \mathrm{diag}(1, p, \ldots, p^{n-1}).
\]
From Lemma~\ref{lemm_paso_0}, the matrix $H_0 = Y_A (\Lambda_p(Y_A)(z^p))^{-1}$ belongs to $\mathrm{GL}_n(E_{0,K})$. Thus, $\widetilde{H}_1 \in \mathrm{GL}_n(E_{0,K})$. Invoking Lemma~\ref{lemm_paso_0} again, we have $\|H_0\| = 1$, which implies $\|\widetilde{H}_1\| \le 1$. Since $\widetilde{H}_1(0) = \mathrm{diag}(1, p, \ldots, p^{n-1})$, we obtain $\|\widetilde{H}_1\| = 1$.

We now show that $\delta \widetilde{H}_1 = A \widetilde{H}_1 - p \widetilde{H}_1 A_1(z^p)$. Recall from statement (b) that $\delta T = A_1 T$, and hence $\delta(T(z^p)) = p A_1(z^p) T(z^p)$. That is, $T(z^p)$ is a fundamental matrix of solutions to $\delta X = p A_1(z^p) X$. Furthermore, statement (b) yields
\[
T(z^p) = \mathrm{diag}\left(1, \frac{1}{p}, \ldots, \frac{1}{p^{n-1}}\right) \Lambda_p(Y_A)(z^p) \, \mathrm{diag}(1, p, \ldots, p^{n-1}) \, z^{p A(0)}.
\]
It then follows directly that
\[
\widetilde{H}_1 T(z^p) = Y_A(z) \, \mathrm{diag}(1, p, \ldots, p^{n-1}) \, z^{p A(0)}.
\]
Since $\mathrm{diag}(1, p, \ldots, p^{n-1}) \, z^{p A(0)} = z^{A(0)} \, \mathrm{diag}(1, p, \ldots, p^{n-1})$, we obtain
\[
\widetilde{H}_1 T(z^p) = Y_A z^{A(0)} \, \mathrm{diag}(1, p, \ldots, p^{n-1}).
\]
Because $Y_A z^{A(0)}$ is a fundamental matrix of solutions to $\delta X = A X$, the matrix $\widetilde{H}_1 T(z^p)$ is also a fundamental matrix of $\delta X = A X$. Consequently,
\begin{align*}
A \widetilde{H}_1 T(z^p) &= \delta(\widetilde{H}_1 T(z^p)) \\
&= (\delta \widetilde{H}_1) T(z^p) + \widetilde{H}_1 \delta(T(z^p)) \\
&= (\delta \widetilde{H}_1) T(z^p) + \widetilde{H}_1 (p A_1(z^p)) T(z^p).
\end{align*}
Canceling the invertible matrix $T(z^p)$ gives $A \widetilde{H}_1 = \delta \widetilde{H}_1 + p \widetilde{H}_1 A_1(z^p)$. Rearranging terms, we conclude
\[
\delta \widetilde{H}_1 = A \widetilde{H}_1 - p \widetilde{H}_1 A_1(z^p). 
\]
\end{proof}

Finally, Steps~1 and~2 follow from the following corollary.

\begin{coro}\label{prop_ant}
Let the assumptions be as in Lemma~\ref{lemm_paso_0}, and let $f(z) \in 1 + z K[[z]]$ be a solution to $\mathcal{L}$. Then, for every integer $m > 0$, there exists a MOM differential operator $\mathcal{L}_m \in E_{0,K}[\delta]$ of order $n$ such that:
\begin{enumerate}[label=(\alph*)]
  \item $\Lambda^m_p(f)$ is a solution to $\mathcal{L}_m$.

  \item If $A_m$ denotes the companion matrix of $\mathcal{L}_m$, then a fundamental matrix of solutions to $\delta X = A_m X$ is given by
  \[
  \mathrm{diag}\left(1, \frac{1}{p^m}, \ldots, \frac{1}{p^{m(n-1)}}\right) \Lambda^m_p(Y_A) \, \mathrm{diag}\left(1, p^m, \ldots, p^{m(n-1)}\right) z^{A(0)}.
  \]

  \item $A_m$ belongs to $\mathcal{M}_K$ and satisfies $\|A_m\| \le 1$.

  \item The matrices $A$ and $p^m A_m(z^{p^m})$ are $E_K$-equivalent. Moreover,
  \[
  \delta \widetilde{H}_m = A \widetilde{H}_m - p^m \widetilde{H}_m A_m(z^{p^m}),
  \]
  where
  \[
  \widetilde{H}_m := Y_A (\Lambda^m_p(Y_A)(z^{p^m}))^{-1} \mathrm{diag}\left(1, p^m, \ldots, p^{m(n-1)}\right) \in \mathrm{GL}_n(E_K) \cap M_n(E_{0,K}),
  \]
  and $\|\widetilde{H}_m\| = 1$.
\end{enumerate}
\end{coro}

\begin{proof}
We proceed by induction on $m \in \mathbb{Z}_{>0}$. 

For $m=1$, we apply Proposition~\ref{lemm_paso_base}. Hence, there exists a MOM differential operator $\mathcal{L}_1 \in E_{0,K}[\delta]$ of order $n$ such that conditions (a)–(d) are satisfied. 

Now, suppose that conditions (a)–(d) hold for some integer $m > 0$. We show that conditions (a)–(d) also hold for $m+1$. By the induction hypothesis, there exists a MOM differential operator $\mathcal{L}_m \in E_{0,K}[\delta]$ of order $n$ such that $\Lambda^m_p(f)$ is a solution to $\mathcal{L}_m$. Moreover, its companion matrix $A_m$ belongs to $\mathcal{M}_K$ with $\|A_m\| \le 1$, and a fundamental matrix of solutions to the system $\delta X = A_m X$ is given by $Y_{A_m} z^{A(0)}$, where
\[
Y_{A_m} = \mathrm{diag}\left(1, \frac{1}{p^m}, \ldots, \frac{1}{p^{m(n-1)}}\right) \Lambda^m_p(Y_A) \, \mathrm{diag}\left(1, p^m, \ldots, p^{m(n-1)}\right).
\]
We also know by the induction hypothesis that the matrix
\[
\widetilde{H}_m = Y_A \bigl(\Lambda^m_p(Y_A)(z^{p^m})\bigr)^{-1} \mathrm{diag}\left(1, p^m, \ldots, p^{m(n-1)}\right)
\]
belongs to $\mathrm{GL}_n(E_K) \cap M_n(E_{0,K})$, satisfies $\|\widetilde{H}_m\| = 1$, and obeys
\begin{equation}\label{eq_equiv}
\delta \widetilde{H}_m = A \widetilde{H}_m - p^m \widetilde{H}_m A_m(z^{p^m}).
\end{equation}

Next, applying Proposition~\ref{lemm_paso_base} to the differential operator $\mathcal{L}_m$, it follows that there exists a MOM differential operator $\mathcal{L}_{m+1} \in E_{0,K}[\delta]$ of order $n$ such that $\Lambda^{m+1}_p(f) = \Lambda_p(\Lambda^m_p(f))$ is a solution to $\mathcal{L}_{m+1}$. The companion matrix $A_{m+1}$ belongs to $\mathcal{M}_K$ with $\|A_{m+1}\| \le 1$, and a fundamental matrix of solutions to $\delta X = A_{m+1} X$ is given by
\[
\mathrm{diag}\left(1, \frac{1}{p}, \ldots, \frac{1}{p^{n-1}}\right) \Lambda_p(Y_{A_m}) \, \mathrm{diag}(1, p, \ldots, p^{n-1}) \, z^{A(0)}.
\]
Substituting the explicit expression for $Y_{A_m}$ into the above formula reveals that a fundamental matrix of solutions to $\delta X = A_{m+1} X$ is given by
\[
\mathrm{diag}\left(1, \frac{1}{p^{m+1}}, \ldots, \frac{1}{p^{(m+1)(n-1)}}\right) \Lambda^{m+1}_p(Y_A) \, \mathrm{diag}\left(1, p^{m+1}, \ldots, p^{(m+1)(n-1)}\right) z^{A(0)}.
\]

Invoking Proposition~\ref{lemm_paso_base} once more, we obtain that the matrix
\[
\widetilde{G}_1 = Y_{A_m} \bigl(\Lambda_p(Y_{A_m})(z^p)\bigr)^{-1} \mathrm{diag}(1, p, \ldots, p^{n-1})
\]
belongs to $\mathrm{GL}_n(E_K) \cap M_n(E_{0,K})$, satisfies $\|\widetilde{G}_1\| = 1$, and fulfills $\delta \widetilde{G}_1 = A_m \widetilde{G}_1 - p \widetilde{G}_1 A_{m+1}(z^p)$.

Evaluating this relation at $z^{p^m}$ and applying the chain rule $\delta(f(z^{p^m})) = p^m (\delta f)(z^{p^m})$ yields
\begin{equation}\label{eq_equiv_nivel_m}
\delta\bigl(\widetilde{G}_1(z^{p^m})\bigr) = p^m A_m(z^{p^m}) \widetilde{G}_1(z^{p^m}) - p^{m+1} \widetilde{G}_1(z^{p^m}) A_{m+1}(z^{p^{m+1}}).
\end{equation}

Setting $\widetilde{H}_{m+1} := \widetilde{H}_m \widetilde{G}_1(z^{p^m})$, we observe that
\[
\widetilde{H}_{m+1} = Y_A \bigl(\Lambda^{m+1}_p(Y_A)(z^{p^{m+1}})\bigr)^{-1} \mathrm{diag}\left(1, p^{m+1}, \ldots, p^{(m+1)(n-1)}\right),
\]
which belongs to $\mathrm{GL}_n(E_K) \cap M_n(E_{0,K})$.

Furthermore, combining Equations~\eqref{eq_equiv} and \eqref{eq_equiv_nivel_m} yields
\[
\delta \widetilde{H}_{m+1} = A \widetilde{H}_{m+1} - p^{m+1} \widetilde{H}_{m+1} A_{m+1}(z^{p^{m+1}}).
\]
Finally, since $\|\widetilde{G}_1\| = 1 = \|\widetilde{H}_m\|$ and $\widetilde{H}_{m+1}(0) = \mathrm{diag}(1, p^{m+1}, \ldots, p^{(m+1)(n-1)})$, we conclude that $\|\widetilde{H}_{m+1}\| = 1$.

Therefore, statements (a)–(d) hold for $m+1$, which completes the induction.
\end{proof}

Before proving assertion (i) of Theorem~\ref{theo_analytic_element}, we make the following remark.

\begin{rema}\label{rem_frob_ac}
Let $\mathcal{L} \in E_{0,K}[\delta]$ be a differential operator of order $n$, let $A$ be the companion matrix of $\mathcal{L}$, and let $H = (h_{i,j})_{1 \le i,j \le n} \in M_n(E_{0,K}[[z]])$ be a matrix such that $\delta H = A H - H B$ for some $B \in M_n(E_{0,K})$. If the column vector $\vec{y} = (y_1, \ldots, y_n)^{\mathrm{T}} \in K[[z]]^n$ is a solution to the system $\delta \vec{y} = B \vec{y}$, then 
\[
h_{1,1} y_1 + h_{1,2} y_2 + \cdots + h_{1,n} y_n
\]
is a solution to $\mathcal{L}$. Indeed, the relation $\delta H = A H - H B$ implies that $H \vec{y}$ is a solution to the system $\delta \vec{x} = A \vec{x}$. Since $A$ is the companion matrix of $\mathcal{L}$, its first entry $h_{1,1} y_1 + h_{1,2} y_2 + \cdots + h_{1,n} y_n$ is a solution to $\mathcal{L}$.
\end{rema}

\subsection{Proof of (i) of Theorem \ref{theo_analytic_element}}
Let us write $f(z) = \sum_{j \ge 0} f_j z^j$. In order to prove that $f(z) \in 1 + z \mathcal{O}_K[[z]]$, it suffices to show that $f_0, f_1, \ldots, f_{p^m-1} \in \mathcal{O}_K$ for all integers $m > 0$. Since $f(z) \in \mathcal{M}\mathcal{F}(K)$, there exists a MOM differential operator $\mathcal{L} \in E_{0,K}[\delta]$ such that $\mathcal{L}$ has a strong Frobenius structure and the coefficients of $\mathcal{L}$ have Gauss norm at most $1$. Let $n$ be the order of $\mathcal{L}$, and let $A$ be its companion matrix. Because $\mathcal{L}$ has a strong Frobenius structure, it follows from Propositions~4.1.2, 4.6.4, and 4.7.2 of \cite{C83} that there exists an invertible matrix $U$ with entries defined on the generic disc $D(t,1^-)$ such that $\delta U = A U$. Furthermore, all eigenvalues of $A(0)$ are zero because $\mathcal{L}$ is MOM at zero. Hence, $A \in \mathcal{M}_K$. Moreover, $\|A\| \le 1$ since the coefficients of $\mathcal{L}$ have Gauss norm at most $1$.

Let $m \ge 1$ be an integer. Then, by Corollary~\ref{prop_ant}, there exists a MOM differential operator $\mathcal{L}_m \in E_{0,K}[\delta]$ of order $n$ such that $\Lambda^m_p(f)$ is a solution to $\mathcal{L}_m$, and its companion matrix $A_m$ belongs to $\mathcal{M}_K$. Furthermore, Corollary~\ref{prop_ant} guarantees that $\|A_m\| \le 1$ and 
\begin{equation}\label{eq_equiv_m}
\delta \widetilde{H}_m = A \widetilde{H}_m - p^m \widetilde{H}_m A_m(z^{p^m}),
\end{equation}
with $\widetilde{H}_m \in \mathrm{GL}_n(E_K) \cap M_n(E_{0,K})$ and $\|\widetilde{H}_m\| = 1$.

Since $A_m$ is the companion matrix of $\mathcal{L}_m$, the column vector
\[
\bigl(\Lambda^m_p(f), \delta(\Lambda^m_p(f)), \ldots, \delta^{n-1}(\Lambda^m_p(f))\bigr)^{\mathrm{T}}
\]
is a solution to the system $\delta \vec{y} = A_m \vec{y}$. Consequently, the vector
\[
\bigl(\Lambda^m_p(f)(z^{p^m}), (\delta(\Lambda^m_p(f)))(z^{p^m}), \ldots, (\delta^{n-1}(\Lambda^m_p(f)))(z^{p^m})\bigr)^{\mathrm{T}}
\]
is a solution to $\delta \vec{y} = p^m A_m(z^{p^m}) \vec{y}$.

Write $\widetilde{H}_m = (\widetilde{h}^{(m)}_{i,j})_{1 \le i,j \le n}$. By Equation~\eqref{eq_equiv_m}, we apply Remark~\ref{rem_frob_ac} to deduce that 
\[
\widetilde{h}^{(m)}_{1,1} \Lambda^m_p(f)(z^{p^m}) + \widetilde{h}^{(m)}_{1,2} (\delta(\Lambda^m_p(f)))(z^{p^m}) + \cdots + \widetilde{h}^{(m)}_{1,n} (\delta^{n-1}(\Lambda^m_p(f)))(z^{p^m})
\]
is a solution to $\mathcal{L}$. According to Remark~\ref{rem_sol}, there exists $c \in \mathbb{C}_p$ such that
\begin{equation}\label{eq_eigenvalue}
\widetilde{h}^{(m)}_{1,1} \Lambda^m_p(f)(z^{p^m}) + \widetilde{h}^{(m)}_{1,2} (\delta(\Lambda^m_p(f)))(z^{p^m}) + \cdots + \widetilde{h}^{(m)}_{1,n} (\delta^{n-1}(\Lambda^m_p(f)))(z^{p^m}) = c f(z).
\end{equation}

Since $f(0) = 1$ and $(\delta^j(\Lambda^m_p(f)))(z^{p^m})(0) = 0$ for all integers $j > 0$, evaluating Equation~\eqref{eq_eigenvalue} at $z = 0$ yields $\widetilde{h}^{(m)}_{1,1}(0) = c$. Moreover, $\widetilde{h}^{(m)}_{1,1}(0) = 1$ because $\widetilde{H}_m(0) = \mathrm{diag}(1, p^m, \ldots, p^{m(n-1)})$, which implies $c = 1$. Note also that $\widetilde{h}_{1,j}^{(m)}(0) = 0$ for all $j \in \{2, \ldots, n\}$. Since $(\delta^j(\Lambda^m_p(f)))(z^{p^m}) \in z^{p^m} K[[z]]$ for all $j > 0$, Equation~\eqref{eq_eigenvalue} reduces to
\begin{equation}\label{eq_reduction}
\widetilde{h}^{(m)}_{1,1} \Lambda^m_p(f)(z^{p^m}) \equiv f(z) \pmod{z^{p^m+1} K[[z]]}.
\end{equation}

Writing $\widetilde{h}^{(m)}_{1,1} = \sum_{j \ge 0} h_j z^j$, Equation~\eqref{eq_reduction} implies that $h_j = f_j$ for every $j \in \{1, \ldots, p^m - 1\}$. Because $\|\widetilde{H}_m\| = 1$, we have $|\widetilde{h}^{(m)}_{1,1}|_{\mathcal{G}} \le 1$, and thus $|f_j| \le 1$ for all $j \in \{1, \ldots, p^m - 1\}$. We conclude that $f_j \in \mathcal{O}_K$ for every integer $j \in \{1, \ldots, p^m - 1\}$. Since $m \ge 1$ was arbitrary, this completes the proof.
$\hfill\square$

\section{Constructing analytic elements and congruences modulo $p^m$}\label{sec_congruences}

The aim of this section is to prove parts (ii) and (iii) of Theorem~\ref{theo_analytic_element}. The cornerstone of our proof is Lemma~\ref{lemm_ecuacion_nivel_1}, which establishes that $\frac{f(z)}{\Lambda_p(f)((z^p))}$ belongs to $E_{0,K}$. As mentioned in Section~\ref{sec_stra}, the proof of this lemma relies on Lemma~\ref{lemm_congruences}, which is very technical but allows us to obtain congruences modulo $p^m$ for any integer $m\geq1$.

Recall that if $K$ is a finite extension of $\mathbb{Q}_p$, then the maximal ideal $\mathfrak{m}_{\mathcal{O}_K}$ of $\mathcal{O}_K$ is a principal ideal. An element $\pi \in K$ is called a \emph{uniformizer} of $ \mathcal{O}_K $ if $\mathfrak{m}_{\mathcal{O}_K} = (\pi)$. We also recall that $K_0(z)$ denotes the ring of rational functions $A(z)/B(z)$, where $A(z), B(z) \in K[z]$ and $B(x) \neq 0$ for all $x \in D_0$.

\subsection{Cartier operator and congruences modulo $p^m$.}

The main result of this section is Lemma~\ref{lemm_ecuacion_nivel_1}. The proof of this lemma relies on some results concerning the Cartier operator and congruences modulo $p^m$.

\begin{lemm}\label{lemm_h_m}
Let the assumptions be as in Lemma~\ref{lemm_paso_0}. For every integer $l \ge 1$, we define 
\[
H_l := \mathrm{diag}\left(1, \frac{1}{p^l}, \ldots, \frac{1}{p^{l(n-1)}}\right) \Lambda^l_p(Y_A) \bigl(\Lambda^{l+1}_p(Y_A)(z^p)\bigr)^{-1} \mathrm{diag}\left(1, p^l, \ldots, p^{l(n-1)}\right).
\] 
Then, $H_l \in \mathrm{GL}_n(E_{0,K})$ and satisfies $\|H_l\| = 1 = \|H_l^{-1}\|$. Moreover, $H_l(0) = I = H_l^{-1}(0)$, and the $(1,1)$-entry of both $H_l$ and $H_l^{-1}$ has Gauss norm $1$.
\end{lemm}

\begin{proof}
Let $l \ge 1$ be an integer. From Corollary~\ref{prop_ant}, it follows that there exists a MOM differential operator $\mathcal{L}_l \in E_{0,K}[\delta]$ of order $n$ such that $\Lambda^l_p(f)$ is a solution to $\mathcal{L}_l$. Let $A_l$ be the companion matrix of $\mathcal{L}_l$. Again by Corollary~\ref{prop_ant}, we know that $A_l$ belongs to $\mathcal{M}_K$, $\|A_l\| \le 1$, and a fundamental matrix of solutions to the system $\delta X = A_l X$ is given by $Y_{A_l} z^{A(0)}$, where 
\[
Y_{A_l} = \mathrm{diag}\left(1, \frac{1}{p^l}, \ldots, \frac{1}{p^{l(n-1)}}\right) \Lambda^l_p(Y_A) \, \mathrm{diag}\left(1, p^l, \ldots, p^{l(n-1)}\right).
\] 
Now, applying Lemma~\ref{lemm_paso_0} to $A_l$, the matrix $H_l = Y_{A_l} \bigl(\Lambda_p(Y_{A_l})(z^p)\bigr)^{-1}$ belongs to $\mathrm{GL}_n(E_{0,K})$ and satisfies $\|H_l\| = 1 = \|H_l^{-1}\|$. Since $H_l(0) = I = H_l^{-1}(0)$, the $(1,1)$-entry of both $H_l$ and $H_l^{-1}$ has Gauss norm $1$.
\end{proof}

The following technical lemma provides congruences modulo $p^m$.

\begin{lemm}\label{lemm_congruences}
Let $l \ge 1$ be an integer, let $H_l$ be the matrix defined in Lemma~\ref{lemm_h_m}, and let $m \ge 2$ be an integer. We define 
\[
\bigl(s_{\mu,\nu}\bigr)_{1 \le \mu, \nu \le n} := \prod_{l=1}^{m-1} \left[ H_l(z^{p^l}) \, \mathrm{diag}\left(1, p, \ldots, p^{n-1}\right) \right].
\]
Then:
\begin{enumerate}[label=(\roman*)]
  \item $|s_{1,1}|_{\mathcal{G}} = 1 = |s_{1,1}(0)|$ and $s_{1,1}$ belongs to $1 + z \mathcal{O}_K[[z^p]]$.

  \item For all $\mu, \nu \in \{2, \ldots, n\}$, we have 
  \[
  s_{\mu,1} s_{1,\nu} - s_{1,1} s_{\mu,\nu} \in p^{m-1} \mathcal{O}_K[[z]].
  \]
\end{enumerate}
\end{lemm}
\begin{proof}
Notice that $s_{\mu,\nu}$ belongs to $E_{0,K}\cap\mathcal{O}_K[[z^p]]$ because, by Lemma~\ref{lemm_h_m}, we deduce that, for any integer $l\geq1$, $||H_l(z^{p^l})||=1$ and $H_l(z^{p^l})\in GL_n(E_{0,K})$. 

For every integer $l\geq1$, we write $H_l(z^{p^l})=(r^{(l)}_{i,j})_{1\leq i,j\leq n}$. Thus, 
\begin{equation}\label{eq_s}
s_{\mu,\nu}=\sum_{j_1=1}^{n}\sum_{j_2=1}^{n}\cdots\sum_{j_{m-2}=1}^n p^{j_1+\cdots+j_{m-2}+\nu-m+1}r^{(1)}_{\mu,j_{m-2}}r^{(2)}_{j_{m-2},j_{m-3}}r^{(3)}_{j_{m-3},j_{m-4}}\cdots r^{(m-1)}_{j_1,\nu}.
\end{equation}
(i).  It is clear that 
\begin{align*}
&s_{1,1}=\\
&r^{(1)}_{1,1}r^{(2)}_{1,1}r^{(3)}_{1,1}\cdots r^{(m-1)}_{1,1}+\sum_{\begin{subarray}{l}1\leq j_1,\ldots, j_{m-2}\leq n,\\
(j_1,\ldots, j_{m-2})\neq(1,\ldots,1)\end{subarray}}p^{j_1+\cdots+j_{m-2}+1-m+1}r^{(1)}_{1,j_{m-2}}r^{(2)}_{j_{m-2},j_{m-3}}\cdots r^{(m-1)}_{j_1,1}.
\end{align*}
 Let $\textbf{J}=(j_1,\ldots, j_{m-2})$ be in $\{1,\ldots, n\}^{m-2}$ such that $\textbf{J}\neq(1,\ldots,1)$. Therefore, there exists $\kappa\in\{1,\ldots, m-2\}$ such that $j_{\kappa}\geq2$. So that, $$j_1+\cdots+j_{m-2}+1-m+1=j_1+\cdots+j_{m-2}-m+2\geq m-3+2-m+2=1>0.$$
From Lemma~\ref{lemm_h_m}, we know that $||H_l||=1$. So, for all $i,j\in\{1,\ldots,n\}$, we have $|r^{(l)}_{i,j}|\leq1$. Consequently, $$\left|\sum_{\begin{subarray}{l}1\leq j_1,\ldots j_{m-2}\leq n,\\
(j_1,\ldots, j_{m-2})\neq(1,\ldots,1)\end{subarray}}p^{j_1+\cdots+j_{m-2}+1-m+1}r^{(1)}_{1,j_{m-2}}r^{(2)}_{j_{m-2},j_{m-3}}\cdots r^{(m-1)}_{j_1,1}\right|_{\mathcal{G}}<1.$$
Furthermore, according to Lemma~\ref{lemm_h_m},  we conclude that, for any $l\geq1$, $|r^{(l)}_{1,1}|=1$. Therefore, we have $|s_{1,1}|_{\mathcal{G}}=1$. Finally, we also have $H_{l}(0)=I$. Hence, if $(j_1,\ldots, j_{m-2})\neq(1,\ldots, 1)$ then $r^{(1)}_{1,j_{m-2}}r^{(2)}_{j_{m-2},j_{m-3}}\cdots r^{(m-1)}_{j_1,1}(0)=0$ and it is clear that $r^{(1)}_{1,1}r^{(2)}_{1,1}r^{(3)}_{1,1}\cdots r^{(m-1)}_{1,1}(0)=1$. Thus, $|s_{1,1}(0)|=1$. So, we have $|s_{1,1}|_{\mathcal{G}}=1=|s_{1,1}(0)|$. Consequently, we obtain $\frac{1}{s_{1,1}}$ belongs to $1+\mathcal{O}_K[[z]]\cap E_{K}$. Even, we have $\frac{1}{s_{1,1}}\in 1+\mathcal{O}_K[[z^p]]$ because, by construction, it is clear that $s_{1,1}\in\mathcal{O}_K[[z^p]]$ and we have already seen that $s_{1,1}(0)=1$.

(ii)  From Equation~\eqref{eq_s}, we have
\begin{small}
\begin{align*}
&s_{\mu,1}s_{1,\nu}-s_{1,1}s_{\mu,\nu}=\\
&\sum_{\begin{subarray}{l}1\leq j_1,\ldots, j_{m-2}\leq n,\\
1\leq t_1,\ldots, t_{m-2}\leq n\end{subarray}}p^{\displaystyle\sum_{\gamma=1}^{m-2} (j_{\gamma}+ t_{\gamma})+\nu-2m+3}r^{(1)}_{\mu,j_{m-2}}r^{(2)}_{j_{m-2},j_{m-3}}\cdots r^{(m-1)}_{j_1,1}r^{(1)}_{1,t_{m-2}}r^{(2)}_{t_{m-2},t_{m-3}}\cdots r^{(m-1)}_{t_1,\nu}\\
&-\\
&\sum_{\begin{subarray}{l}1\leq j_1,\ldots, j_{m-2}\leq n,\\
1\leq t_1,\ldots, t_{m-2}\leq n\end{subarray}}p^{\displaystyle\sum_{\gamma=1}^{m-2} (j_{\gamma}+ t_{\gamma})+\nu-2m+3}r^{(1)}_{1,j_{m-2}}r^{(2)}_{j_{m-2},j_{m-3}}\cdots r^{(m-1)}_{j_1,1}r^{(1)}_{\mu,t_{m-2}}r^{(2)}_{t_{m-2},t_{m-3}}\cdots r^{(m-1)}_{t_1,\nu}.
\end{align*}
\end{small}
Let $\textbf{J}=(j_1,\ldots, j_{m-2})$ and $\textbf{T}=(t_1,\ldots, t_{m-2})$ be in $\{1,\ldots,n\}^{m-2}$ such that $j_{\kappa}=t_{\kappa}$ for some  $\kappa\in\{1,\ldots,m-2\}$.  Let us consider $\textbf{J'}=(j'_1,\ldots, j'_{m-2})$ and $\textbf{T'}=(t'_1,\ldots,t'_{m-2})$ defined as follows $$(j'_1,\ldots, j'_{\kappa-1}, j'_{\kappa}, j'_{\kappa+1}, j'_{\kappa+2}\ldots, j'_{m-2})=(j_1,\ldots, j_{\kappa-1},j_{\kappa},t_{\kappa+1}, t_{\kappa+2}\ldots, t_{m-2})$$ and $$(t'_1,\ldots, t'_{\kappa-1}, t'_{\kappa},t'_{\kappa+1}, t'_{\kappa+2},\ldots,t'_{m-2})=(t_1,\ldots, t_{\kappa-1},t_{\kappa},j_{\kappa+1},\ldots, j_{m-2}.)$$
So that, we have
\begin{small}
\begin{align*}
&\bigg[r^{(1)}_{\mu,j_{m-2}}r^{(2)}_{j_{m-2},j_{m-3}}\cdots r^{(m-1)}_{j_1,1}\bigg]\bigg[r^{(1)}_{1,t_{m-2}}r^{(2)}_{t_{m-2},t_{m-3}}\cdots r^{(m-1)}_{t_1,\nu}\bigg]=\\
&\bigg[r^{(1)}_{\mu,j_{m-2}}\cdots r^{(m-\kappa-1)}_{j_{\kappa+1},j_{\kappa}}r^{(m-\kappa)}_{j_{\kappa},j_{\kappa-1}}r^{(m-\kappa+1)}_{j_{\kappa-1},j_{\kappa-2}}\cdots r^{(m-1)}_{j_1,1}\bigg]\bigg[r^{(1)}_{1,t_{m-2}}\cdots r^{(m-\kappa-1)}_{t_{\kappa+1},t_{\kappa}}r^{(m-\kappa)}_{t_{\kappa},t_{\kappa-1}}r^{(m-\kappa+1)}_{t_{\kappa-1},t_{\kappa-2}}\cdots r^{(m-1)}_{t_1,\nu}\bigg]=\\
&\bigg[r^{(1)}_{1,t_{m-2}}\cdots r^{(m-\kappa-1)}_{t_{\kappa+1},t_{\kappa}}r^{(m-\kappa)}_{j_{\kappa},j_{\kappa-1}}r^{(m-\kappa+1)}_{j_{\kappa-1},j_{\kappa-2}}\cdots r^{(m-1)}_{j_1,1}\bigg]\bigg[r^{(1)}_{\mu,j_{m-2}}\cdots r^{(m-\kappa-1)}_{j_{\kappa+1},j_{\kappa}}r^{(m-\kappa)}_{t_{\kappa},t_{\kappa-1}}r^{(m-\kappa+1)}_{t_{\kappa-1},t_{\kappa-2}}\cdots r^{(m-1)}_{t_1,\nu}\bigg]=\\
&\bigg[r^{(1)}_{1,j'_{m-2}}\cdots r^{(m-\kappa-1)}_{j'_{\kappa+1},j'_{\kappa}}r^{(m-\kappa)}_{j'_{\kappa},j'_{\kappa-1}}r^{(m-\kappa+1)}_{j'_{\kappa-1},j'_{\kappa-2}}\cdots r^{(m-1)}_{j'_1,1}\bigg]\bigg[r^{(1)}_{\mu,t'_{m-2}}\cdots r^{(m-\kappa-1)}_{t'_{\kappa+1},t'_{\kappa}}r^{(m-\kappa)}_{t'_{\kappa},t'_{\kappa-1}}r^{(m-\kappa+1)}_{t'_{\kappa-1},t'_{\kappa-2}}\cdots r^{(m-1)}_{t'_1,\nu}\bigg].
\end{align*}
\end{small}
Moreover, we also have $\displaystyle\sum_{\gamma=1}^{m-2} (j_{\gamma}+ t_{\gamma})=\displaystyle\sum_{\gamma=1}^{m-2} (j'_{\gamma}+ t'_{\gamma})$. Thus, we obtain
\begin{align*}
&p^{\sum (j_{\gamma}+ t_{\gamma})+\nu-2m+3}r^{(1)}_{\mu,j_{m-2}}r^{(2)}_{j_{m-2},j_{m-3}}\cdots r^{(m-1)}_{j_1,1}r^{(1)}_{1,t_{m-2}}r^{(2)}_{t_{m-2},t_{m-3}}\cdots r^{(m-1)}_{t_1,\nu}-\\
&p^{\sum (j'_{\gamma}+ t'_{\gamma})+\nu-2m+3}r^{(1)}_{1,j'_{m-2}}r^{(2)}_{j'_{m-2},j'_{m-3}}\cdots r^{(m-1)}_{j'_1,1}r^{(1)}_{\mu,t'_{m-2}}r^{(2)}_{t'_{m-2},t'_{m-3}}\cdots r^{(m-1)}_{t'_1,\nu}\\
&=0\\
%&p^{\sum (j_{\gamma}+ t_{\gamma})+\nu-2m+3}\bigg[r^{(1)}_{\mu,j_{m-2}}r^{(2)}_{j_{m-2},j_{m-3}}\cdots r^{(m-1)}_{j_1,1}r^{(1)}_{1,t_{m-2}}r^{(2)}_{t_{m-2},t_{m-3}}\cdots r^{(m-1)}_{t_1,\nu}-\\
%&r^{(1)}_{\mu,j'_{m-2}}r^{(2)}_{j'_{m-2},j'_{m-3}}\cdots r^{(m-1)}_{j'_1,1}r^{(1)}_{1,t'_{m-2}}r^{(2)}_{t'_{m-2},t'_{m-3}}\cdots r^{(m-1)}_{t'_1,\nu}\bigg]\\
%&=0
\end{align*}
Consequently, 
\begin{small}
\begin{align*}
&s_{\mu,1}s_{1,\nu}-s_{1,1}s_{\mu,\nu}=\\
&\sum_{\substack{1\leq j_1,\ldots, j_{m-2}\leq n\\
1\leq t_1,\ldots, t_{m-2}\leq n\\
j_i\neq t_i\forall i}}p^{\displaystyle\sum_{\gamma=1}^{m-2} (j_{\gamma}+ t_{\gamma})+\nu-2m+3}r^{(1)}_{\mu,j_{m-2}}r^{(2)}_{j_{m-2},j_{m-3}}\cdots r^{(m-1)}_{j_1,1}r^{(1)}_{1,t_{m-2}}r^{(2)}_{t_{m-2},t_{m-3}}\cdots r^{(m-1)}_{t_1,\nu}\\
&-\\
&\sum_{\substack{1\leq j_1,\ldots, j_{m-2}\leq n\\
1\leq t_1,\ldots, t_{m-2}\leq n\\
j_i\neq t_i\forall i}}p^{\displaystyle\sum_{\gamma=1}^{m-2} (j_{\gamma}+ t_{\gamma})+\nu-2m+3}r^{(1)}_{1,j_{m-2}}r^{(2)}_{j_{m-2},j_{m-3}}\cdots r^{(m-1)}_{j_1,1}r^{(1)}_{\mu,t_{m-2}}r^{(2)}_{t_{m-2},t_{m-3}}\cdots r^{(m-1)}_{t_1,\nu}.
\end{align*}
\end{small}
Let $\textbf{J}=(j_1,\ldots, j_{m-2})$ and $\textbf{T}=(t_1,\ldots, t_{m-2})$ be in $\{1,\ldots,n\}^{m-2}$ such that, for all $\kappa\in\{1,\ldots,m-2\}$, $j_{\kappa}\neq t_{\kappa}$. Hence, for all $\kappa\in\{1,\ldots, m-2\}$, $j_{\kappa}+t_{\kappa}\geq 3$. Therefore,  $$\sum_{\gamma=1}^{m-2} j_{\gamma}+t_{\gamma}\geq 3m-6.$$
As $\nu\geq2$ then $$\sum_{\gamma=1}^{m-2} (j_{\gamma}+t_{\gamma})+\nu-2m+3\geq m-1.$$
Therefore,
$$s_{\mu,1}s_{1,\nu}-s_{1,1}s_{\mu,\nu}\in p^{m-1}\mathcal{O}_K[[z]].$$
\end{proof}
We are now ready to show that $\frac{f(z)}{\Lambda_p(f)(z^p)}$ belongs to $E_p$.
\begin{lemm}\label{lemm_ecuacion_nivel_1}
Let the assumptions be as in Lemma~\ref{lemm_paso_0}, and let $f(z)\in 1+zK[[z]]$ be a solution of $\mathcal{L}$. Then, for all integers $m>0$, there exists $D_m(z)\in K_0(z)$ with norm 1 such that $$D_m(z)(\Lambda_p(f))(z^p)\equiv f(z)\bmod\pi^m\mathcal{O}_K[[z]],$$
where $\pi$ is a uniformizer of $\mathcal{O}_K$.
\end{lemm}

\begin{proof}
We split the proof into two cases.

\textbullet\quad\emph{Case $m=1$}. According to Proposition~\ref{lemm_paso_base}, there exists a MOM differential operator $\mathcal{L}_1\in E_{0,K}[\delta]$ such that $\Lambda_p(f)$ is a solution of $\mathcal{L}_1$. Let $A_1$ be the companion matrix of $\mathcal{L}_1$ and let $B_0=Y_A(\Lambda(Y_A)(z^p))^{-1}$, where $Y_A$ is the uniform part of $\delta X=AX$. Again by Proposition~\ref{lemm_paso_base}, the matrix $\widetilde{H}_1=B_0 \mathrm{diag}\left(1,p,\ldots, p^{n-1}\right)$ belongs to $GL_n(E_K)\cap M_n(E_{0,K})$ and 
\begin{equation}\label{eq_nivel_1}
\delta \widetilde{H}_1= A\widetilde{H}_1-p\widetilde{H}_1A_1(z^p).
\end{equation}
Since $\Lambda_p(f)$ is a solution of $\mathcal{L}_1$, the vector $((\Lambda_p(f))(z^p),(\delta\Lambda_p(f))(z^p),\ldots, (\delta^{n-1}\Lambda_p(f))(z^p))^t$ is a solution of  $\delta\vec{y}=pA_1(z^p)\vec{y}$.  Let us write $B_0=(b_{i,j})_{1\leq i,\leq n}$. Since $A$ is the companion matrix of $\mathcal{L}$, Remark~\ref{rem_frob_ac} together with Equation~\eqref{eq_nivel_1} implies that $$b_{1,1}(\Lambda_p(f))(z^p)+pb_{1,2}(\delta\Lambda_p(f))(z^p)+\cdots+ p^{n-1}b_{1,n} (\delta^{n-1}\Lambda_p(f))(z^p)$$
is a solution of $\mathcal{L}$.  Furthermore, this solution belongs to $K[[z]]$ because, by Lemma~\ref{lemm_paso_0}, $B_0\in GL_n(E_{0,K})$ and it is clear that $\Lambda_p(f)\in 1+zK[[z]]$.

As $\mathcal{L}$ is MOM at zero, from Remark~\ref{rem_sol}, there exists $c\in K$ such that 
\begin{equation}\label{eq_1}
b_{1,1}(\Lambda_p(f))(z^p)+pb_{1,2}(\delta\Lambda_p(f))(z^p)+\cdots+ p^{n-1}b_{1,n} (\delta^{n-1}\Lambda_p(f))(z^p)=cf(z).
\end{equation}
But $c=1$ because $B_0(0)=I$ and $f(0)=1$. 

Since $B_0\in GL_n(E_{0,K})$ and $p\in\pi\mathcal{O}_K$, we can reduce Equation~\eqref{eq_1} modulo $\pi$ and, we obtain $$b_{1,1}(\Lambda_p(f))(z^p)\equiv f(z)\bmod\pi\mathcal{O}_K[[z]].$$
We have $|b_{1,1}|=1$ because $||B_0||=1$ and $b_{1,1}(0)=1$. Thus, since $b_{1,1}\in E_{0,K}$, there exists $D_1\in K_0(z)$ with norm equal to 1 such that $b_{1,1}\equiv D_1\bmod\pi\mathcal{O}_K[[z]]$. Consequently, $$D_{1}(\Lambda_p(f))(z^p)\equiv f(z)\bmod\pi\mathcal{O}_K[[z]].$$

\textbullet\quad\emph{Case $m\geq2$}. By Corollary~\ref{prop_ant}, there exists a MOM differential operator $\mathcal{L}_m\in E_{0,K}[\delta]$ of order $n$ such that $\Lambda_p^m(f)$ is a solution of $\mathcal{L}_m$. Let $A_m$ be the companion matrix of $\mathcal{L}_m$. Then, by invoking Corollary~\ref{prop_ant} again, $A_m$ belongs to $\mathcal{M}_K$ and $||A_m||\leq1$. Let us consider
\[\widetilde{H}_m=Y_A(\Lambda^m_p(Y_A)(z^{p^m}))^{-1} \mathrm{diag}\left(1,p^m,\ldots, p^{m(n-2)},p^{m(n-1)}\right).\]
By (d) of Corollary~\ref{prop_ant},  $\widetilde{H}_m$ belongs to $ GL_n(E_K)\cap M_n(E_{0,K})$ , $||\widetilde{H}_m||=1$ and 
\begin{equation}\label{eq_nivel_m}
\delta(\widetilde{H}_m)=A\widetilde{H}_m-p^m\widetilde{H}_mA_m(z^{p^m}).
\end{equation}
Let us write $\widetilde{H}_m=(\widetilde{h}^{(m)}_{i,j})_{1\leq i,j\leq n}$. Hence $|\widetilde{h}^{(m)}_{1,1}|_{\mathcal{G}}\leq1$. But, $\widetilde{h}^{(m)}_{1,1}(0)=1$ given that $\widetilde{H}_m(0)= \mathrm{diag}\left(1,p^m,\ldots, p^{m(n-1)}\right)$. Thus, $|\widetilde{h}^{(m)}_{1,1}|_{\mathcal{G}}=1$. 

We now prove this case in several steps.

\textbullet\quad\textit{Step 1}. We have $$\widetilde{h}^{(m)}_{1,1}\Lambda^m_p(f)(z^{p^m})+\widetilde{h}^{(m)}_{1,2}\delta(\Lambda^m_p(f))(z^{p^m})+\cdots+\widetilde{h}^{(m)}_{1,n}\delta^{n-1}(\Lambda^m_p(f))(z^{p^m})=f(z).$$
Indeed, as $\Lambda^m_p(f)$ is a solution of $\mathcal{L}_m$ and $A_m$ is the companion matrix of $\mathcal{L}_m$ then the vector given by $(\Lambda^m_p(f),\delta(\Lambda^m_p(f)),\ldots,\delta^{n-1}(\Lambda^m_p(f)))^t$ is a solution of $\delta\vec{y}=A_m\vec{y}$. Thus, it follows that the vector $(\Lambda^m_p(f)(z^{p^m}),\delta(\Lambda^m_p(f))(z^{p^m}),\ldots,\delta^{n-1}(\Lambda^m_p(f))(z^{p^m}))^t$ is a solution of $\delta\vec{y}=p^mA_m(z^{p^m})\vec{y}$. Therefore, thanks to Equation~\eqref{eq_nivel_m} we can apply Remark~\ref{rem_frob_ac} and we deduce that $$\widetilde{h}^{(m)}_{1,1}\Lambda^m_p(f)(z^{p^m})+\widetilde{h}^{(m)}_{1,2}\delta(\Lambda^m_p(f))(z^{p^m})+\cdots+\widetilde{h}^{(m)}_{1,n}\delta^{n-1}(\Lambda^m_p(f))(z^{p^m})$$
is a solution of $\mathcal{L}$. In addition, this solution is in $K[[z]]$ because $\widetilde{H}_m\in M_n(E_{0,K})$ and it is clear that $\Lambda_p(f)\in 1+zK[[z]]$.

Thus, as $\mathcal{L}$ is MOM at zero, according to Remark~\ref{rem_sol}, there exists $r\in K$ such that $$\widetilde{h}^{(m)}_{1,1}\Lambda^m_p(f)(z^{p^m})+\widetilde{h}^{(m)}_{1,2}\delta(\Lambda^m_p(f))(z^{p^m})+\cdots+\widetilde{h}^{(m)}_{1,n}\delta^{n-1}(\Lambda^m_p(f))(z^{p^m})=rf(z).$$
Since $f(0)=1$ and, for all integers $j>0$, $\delta^{j}(\Lambda^m_p(f))(0)=0$, then $\widetilde{h}^{(m)}_{1,1}(0)=r$. But, we know that $\widetilde{h}^{(m)}_{1,1}(0)=1$. Consequently, $r=1$. That finishes the proof of Step 1.
%\begin{equation}\label{eq_valor_propio}
%\widetilde{h}^{(m)}_{1,1}\Lambda^m_p(f)(z^{p^m})+\widetilde{h}^{(m)}_{1,2}\delta(\Lambda^m_p(f))(z^{p^m})+\cdots+\widetilde{h}^{(m)}_{1,n}\delta^{n-1}(\Lambda^m_p(f))(z^{p^m})=f(z).
%\end{equation}

\textbullet\quad\textit{Step 2}. We now proceed to show that, for all $\nu\in\{2,\ldots,n\}$, there exists $g_{\nu}(z)\in E_{0,K}\cap\mathcal{O}_K[[z^p]]$ such that 
%\begin{equation}\label{eq_super_m}
$$\widetilde{h}^{(m)}_{1,\nu}\equiv\widetilde{h}^{(m)}_{1,1}g_{\nu}(z)\bmod\mathfrak{\pi}^m\mathcal{O}_K[[z]].$$
%\end{equation}

For each integer $l\geq0$, we set $B_l=\Lambda^l_p(Y_A)(\Lambda^{l+1}_p(Y_A)(z^p))^{-1}$. Thus, we have $$Y_A(\Lambda^m_p(Y_A)(z^{p^m}))^{-1}=\prod_{l=0}^{m-1}B_l(z^{p^l}).$$
For every integer $l\geq1$, we define $$H_l= \mathrm{diag}\left(1,1/p^l,\ldots, 1/p^{l(n-1)}\right)\cdot\Lambda^l_p(Y_A)(\Lambda^{l+1}_p(Y_A)(z^{p}))^{-1}\cdot  \mathrm{diag}\left(1, p^l,\ldots, p^{l(n-1)}\right).$$
Hence \[B_l(z^{p^l})= \mathrm{diag}\left(1,p^l,\ldots, p^{l(n-1)}\right)H_l(z^{p^l}) \mathrm{diag}\left(1,1/p^l,\ldots, 1/p^{l(n-1)}\right).
\]
Thanks to Lemma~\ref{lemm_h_m}, $H_l\in GL_n(E_{0,K})$, $H_l(0)=I=H_l(0)^{-1}$, and $||H_l||=1=||H_l^{-1}||$. 

In addition, we also have 
\begin{equation}\label{eq_h_m_tilde}
\widetilde{H}_m=B_0\cdot \mathrm{diag}\left(1,p,\ldots, p^{n-1}\right)
\prod_{l=1}^{m-1}\left[H_l(z^{p^l}) \mathrm{diag}\left(1,p,\ldots, p^{n-1}\right)\right].
\end{equation}
Let us write  \[(s_{\mu,\nu})_{1\leq\mu, \nu\leq n}=\prod_{l=1}^{m-1}\left[H_l(z^{p^l}) \mathrm{diag}\left(1,p,\ldots, p^{n-1}\right)\right].\]

We recall that $B_0=(b_{i,j})_{1\leq i,j\leq n}$. Then, from Equation~\eqref{eq_h_m_tilde}, we get that, for all $\nu\in\{1,\ldots,n\}$, $$\widetilde{h}^{(m)}_{1,\nu}=b_{1,1}s_{1,\nu}+pb_{1,2}s_{2,\nu}+\cdots+p^{n-1}b_{1,n}s_{n,\nu}.$$
We now prove that, for all $\nu\in\{2,\ldots,n\}$, 
\begin{equation}\label{eq_modulo_k}
\widetilde{h}^{(m)}_{1,1}\frac{s_{1,\nu}}{s_{1,1}}\equiv\widetilde{h}^{(m)}_{1,\nu}\bmod\mathfrak{\pi}^m\mathcal{O}_K[[z]].
\end{equation}
It is clear that $$\widetilde{h}^{(m)}_{1,1}\frac{s_{1,\nu}}{s_{1,1}}-\widetilde{h}^{(m)}_{1,\nu}=\sum_{\mu=2}^{n}p^{\mu-1}b_{1,\mu}\left(s_{\mu,1}\frac{s_{1,\nu}}{s_{1,1}}-s_{\mu,\nu}\right).$$
Since, by definition, $B_0=Y_A(\Lambda_p(Y_A)(z^p))^{-1}$, it follows from Lemma~\ref{lemm_paso_0} that $B_0\in GL_n(E_{0,K})$ and that $||B_0||=1$. In particular, for all $\mu\in\{2,\ldots,n\}$, $b_{1,\mu}\in\mathcal{O}_K[[z]]$. Further, $s_{\mu,\nu}\in E_{0,K}\cap\mathcal{O}_K[[z]]$ for all $1\leq\mu, \nu\leq n$ because, according to Lemma~\ref{lemm_h_m}, $H_l\in GL_n(E_{0,K})$ and $||H_l||=1$ for any integer $l>0$. Finally, from (i) of Lemma~\ref{lemm_congruences}, we deduce that $1/s_{1,1}\in 1+\mathcal{O}_K[[z]]$. Therefore,  for all $\mu\in\{2,\ldots,n\}$, $$p^{\mu-1}b_{1,\mu}\left(s_{\mu,1}\frac{s_{1,\nu}}{s_{1,1}}-s_{\mu,\nu}\right)\in\mathcal{O}_K[[z]].$$ 

Hence, to prove the equality described by Equation~\eqref{eq_modulo_k} it is sufficient to show that, for all $\mu\in\{2,\ldots,n\}$, $p^{\mu-1}b_{1,\mu}\left(s_{\mu,1}\frac{s_{1,\nu}}{s_{1,1}}-s_{\mu,\nu}\right)\in\pi^m\mathcal{O}_K[[z]]$. As $s_{1,1}$ belongs to $1+z\mathcal{O}_K[[z]]$ then, it follows that $p^{\mu-1}b_{1,\mu}\left(s_{\mu,1}\frac{s_{1,\nu}}{s_{1,1}}-s_{\mu,\nu}\right)\in\pi^m\mathcal{O}_K[[z]]$ if and only if $p^{\mu-1}b_{1,\mu}\left(s_{\mu,1}s_{1,\nu}-s_{1,1}s_{\mu,\nu}\right)\in\pi^m\mathcal{O}_K[[z]]$. 
But, according to (ii) of Lemma~\ref{lemm_congruences}, we know that, for all $\mu\in\{2,\ldots, n\}$ and for all $\nu\in\{2,\ldots, n\}$, $$\left(s_{\mu,1}s_{1,\nu}-s_{1,1}s_{\mu,\nu}\right)\in p^{m-1}\mathcal{O}_K[[z]].$$
Notice that $m-1+\mu-1\geq m$ because $\mu\geq 2$. Hence, $$p^{\mu-1}b_{1,\mu}\left(s_{\mu,1}s_{1,\nu}-s_{1,1}s_{\mu,\nu}\right)\in p^m\mathcal{O}_K[[z]].$$
Since $p\in\pi\mathcal{O}_K$, $$p^{\mu-1}b_{1,\mu}\left(s_{\mu,1}s_{1,\nu}-s_{1,1}s_{\mu,\nu}\right)\in\pi^m\mathcal{O}_K[[z]].$$ Therefore, for all $\nu\in\{2,\ldots,n\}$,  \begin{equation*}
\widetilde{h}^{(m)}_{1,1}\frac{s_{1,\nu}}{s_{1,1}}\equiv\widetilde{h}^{(m)}_{1,\nu}\bmod\pi^m\mathcal{O}_K[[z]].
\end{equation*}
We put $g_{\nu}=\frac{s_{1,\nu}}{s_{1,1}}$ for $2\leq \nu\leq n$. Thus, $g_{\nu}$ belongs to $E_{0,K}\cap\mathcal{O}_K[[z^p]]$. That completes the proof fo Step 2.

\textbullet\quad\textit{Step 3} There exists  $D_m\in K_0(z)$ with norm 1 such that $$D_m(\Lambda_p(f))(z^p)\equiv f\bmod\pi^m\mathcal{O}_K[[z]].$$
It follows from Step 1 and 2 that $$\widetilde{h}^{(m)}_{1,1}\left(\Lambda^m_p(f)(z^{p^m})+g_2(z)\delta(\Lambda^m_p(f))(z^{p^m})+\cdots+\delta^{n-1}g_n(z)\Lambda^m_p(f))(z^{p^m})\right)\equiv f\bmod\pi^m\mathcal{O}_K[[z]].$$
But, $\Lambda^m_p(f)(z^{p^m})+g_2(z)\delta(\Lambda^m_p(f))(z^{p^m})+\cdots+g_n(z)\delta^{n-1}(\Lambda^m_p(f))(z^{p^m})\in\mathcal{O}_K[[z^p]]$ because, by construction, for all $\nu\in\{2,\ldots,n\}$, $g_{\nu}\in{O}_K[[z^p]]$. In addition, from (i) of Theorem~\ref{theo_analytic_element}, we know that $f(z)\in 1+z\mathcal{O}_K[[z]]$ and thus, $\delta^{j}(\Lambda^m_p(f))\in\mathcal{O}_K[[z]]$ for all integers $j\geq0$. So, let $d(z)$ be in $\mathcal{O}_K[[z]]$ such that $$d(z^p)=\Lambda^m_p(f)(z^{p^m})+g_2(z)\delta(\Lambda^m_p(f))(z^{p^m})+\cdots+\delta^{n-1}g_n(z)(\Lambda^m_p(f))(z^{p^m}).$$
Therefore, $$\widetilde{h}^{(m)}_{1,1} d(z^p)\equiv f\bmod\pi^m\mathcal{O}_K[[z]].$$
By applying $\Lambda_p$ to the previous equality, we get $$\Lambda_p\left(\widetilde{h}^{(m)}_{1,1}\right)d(z)\equiv\Lambda_p(f)\bmod\pi^m\mathcal{O}_K[[z]].$$
As $\widetilde{h}^{(m)}_{1,1}\in 1+z\mathcal{O}_K[[z]]$ then $\Lambda_p\left(\widetilde{h}^{(m)}_{1,1}\right)\in 1+z\mathcal{O}_K[[z]]$ and thus $\Lambda_p\left(\widetilde{h}^{(m)}_{1,1}\right)\bmod\pi^m\mathcal{O}_K$ is a unit element of the ring $(\mathcal{O}_K/\pi^m)[[z]]$. So $$d(z^p)\equiv\frac{(\Lambda_p(f))(z^p)}{(\Lambda_p(\widetilde{h}^{(m)}_{1,1}))(z^p)}\bmod\pi^m\mathcal{O}_K[[z]].$$
Consequently, 
$$\frac{\widetilde{h}^{(m)}_{1,1}}{\Lambda_p((\widetilde{h}^{(m)}_{1,1}))(z^p)}(\Lambda_p(f))(z^p)\equiv f\bmod\pi^m\mathcal{O}_K[[z]].$$

We know that $\widetilde{h}^{(m)}_{1,1}\in E_{0,K}$ and, from Proposition~5.1 of \cite{C86}, we deduce that $\Lambda_p(\widetilde{h}^{(m)}_{1,1})\in E_{0,K}$. Since $|\widetilde{h}^{(m)}_{1,1}|_{\mathcal{G}}=1=|\widetilde{h}^{(m)}_{1,1}(0)|$, we have $|\Lambda_p(\widetilde{h}^{(m)}_{1,1})|_{\mathcal{G}}=1=|\Lambda_p(\widetilde{h}^{(m)}_{1,1})(0)|$. Thus, $\widetilde{h}^{(m)}_{1,1}\big/\Lambda_p((\widetilde{h}^{(m)}_{1,1}))(z^p)$ is an element of $E_{0,K}$ with norm 1. For this reason, there exists $D_m\in K_0(z)$ with norm 1 such that $\widetilde{h}^{(m)}_{1,1}\big/\Lambda_p((\widetilde{h}^{(m)}_{1,1}))(z^p)\equiv D_m\bmod\pi^m\mathcal{O}_K[[z]]$. So that, $$D_m(\Lambda_p(f))(z^p)\equiv f\bmod\pi^m\mathcal{O}_K[[z]].$$

\end{proof}

As a consequence of Lemma~\ref{lemm_ecuacion_nivel_1}, we obtain the following result.

\begin{coro}\label{prop_induction}
Let the assumptions be as in Lemma~\ref{lemm_paso_0}, and let $f(z) \in 1 + z K[[z]]$ be a solution of $\mathcal{L}$. Then, for every integer $k \ge 1$, the quotient $\frac{f(z)}{(\Lambda_p^k(f))(z^{p^k})}$ belongs to $E_{0,K}$.
\end{coro}

\begin{proof}
Let $m \ge 1$ be an integer. Corollary~\ref{prop_ant} ensures that $\Lambda_p^m(f)$ is a solution of a MOM differential operator $\mathcal{L}_m \in E_{0,K}[\delta]$ whose companion matrix $A_m$ satisfies $A_m \in \mathcal{M}_K$ and $\|A_m\| \le 1$. Thus, Lemma~\ref{lemm_ecuacion_nivel_1} implies that $\frac{\Lambda_p^m(f(z))}{(\Lambda_p^{m+1}(f))(z^p)}$ belongs to $E_{0,K}$. Substituting $z$ by $z^{p^m}$, we deduce that for every integer $m \ge 1$, the quotient
\[
\frac{(\Lambda_p^m(f))(z^{p^m})}{(\Lambda_p^{m+1}(f))(z^{p^{m+1}})} \in E_{0,K}.
\]
Now, for any integer $k \ge 1$, we have the telescoping factorisation
\[
\frac{f(z)}{(\Lambda_p^k(f))(z^{p^k})} = \frac{f(z)}{(\Lambda_p(f))(z^p)} \cdot \frac{(\Lambda_p(f))(z^p)}{(\Lambda_p^2(f))(z^{p^2})} \cdots \frac{(\Lambda_p^{k-1}(f))(z^{p^{k-1}})}{(\Lambda_p^k(f))(z^{p^k})}.
\]
Since each factor belongs to $E_{0,K}$, we conclude that $\frac{f(z)}{(\Lambda_p^k(f))(z^{p^k})} \in E_{0,K}$.
\end{proof}

An important consequence of this corollary is given by the following proposition, which plays a fundamental role in the proof of item (ii) of Theorem~\ref{theo_analytic_element}.

\begin{prop}\label{prop_rational}
Let $K$ be a Frobenius field and let $\pi$ be a uniformizer of $\mathcal{O}_K$. Let $\mathcal{L} \in E_{0,K}[\delta]$ be a MOM differential operator with companion matrix $A$, and let $f(z) \in 1 + z K[[z]]$ be a solution of $\mathcal{L}$. If $\|A\| \le 1$ and $\mathcal{L}$ admits a strong Frobenius structure of period $h$, then there exists a finite extension $\mathbf{K}$ of $K$ such that for every integer $k \ge 0$, there exists a rational function $Q_k(z) \in \mathbf{K}_0(z)$ with $\|Q_k\|_{\mathcal{G}} = 1$ satisfying
\[
Q_k(z) \equiv \frac{\Lambda_p^{kh}(f)}{f} \bmod{\pi_{\mathbf{K}}^{kh} \mathcal{O}_{\mathbf{K}}[[z]]},
\]
where $\pi_{\mathbf{K}}$ is a uniformizer of $\mathcal{O}_{\mathbf{K}}$.
\end{prop}

\begin{rema}\label{rem_k_0}
Let $K$ be an extension of $\mathbb{Q}_p$ and let $P(z) \in K_0(z)$. If $|P(z)|_{\mathcal{G}} \le 1$, then $P(z) \in \mathcal{O}_K[[z]]$. Indeed, since $0$ is not a pole of $P(z)$, we can write $P(z) = \sum_{n \ge 0} a_n z^n$ with $a_n \in K$ for all $n \ge 0$. Since $|P(z)|_{\mathcal{G}} \le 1$, we have $\sup_{n \ge 0} |a_n| = |P(z)|_{\mathcal{G}} \le 1$. Thus, $|a_n| \le 1$ for all $n \ge 0$, which yields $P(z) \in \mathcal{O}_K[[z]]$.
\end{rema}

\begin{rema}\label{rem_ast}
Let $K$ be a finite extension of $\mathbb{Q}_p$. Then $E_{0,K} = E_{0,p} \cap K[[z]]$. We first show that $\mathbb{C}_{0,p}(z) \cap K[[z]] = K_0(z)$. Recall that $\mathbb{C}_{0,p}(z)$ denotes the set of rational functions $A(z)/B(z) \in \mathbb{C}_p(z)$ such that $B(x) \neq 0$ for all $x \in D_0$. 

Clearly, $K_0(z) \subset \mathbb{C}_{0,p}(z) \cap K[[z]]$. Conversely, let $P(z) = \frac{\sum_i a_i z^i}{\sum_j b_j z^j} \in \mathbb{C}_{0,p}(z)$ be such that $P(z) = \sum_{i \ge 0} r_i z^i \in K[[z]]$. For every $\sigma \in \mathrm{Gal}(\mathbb{C}_p/K)$, we define $P^{\sigma}(z) = \frac{\sum_i \sigma(a_i) z^i}{\sum_j \sigma(b_j) z^j}$. Expanding as a power series, we obtain $P^{\sigma}(z) = \sum_{i \ge 0} \sigma(r_i) z^i$. Since $\sigma|_K$ is the identity, $r_i \in K$ implies $P^{\sigma}(z) = P(z)$. Consequently, by the Ax-Sen-Tate Theorem, the field of elements in $\mathbb{C}_p$ fixed by $\mathrm{Gal}(\mathbb{C}_p/K)$ is $K$, which implies that the polynomials $\sum_i a_i z^i$ and $\sum_j b_j z^j$ can be chosen in $K[z]$. Hence, $P(z) \in K_0(z)$. 

So we have $\mathbb{C}_{0,p}(z) \cap K[[z]] = K_0(z)$, and taking completions with respect to the Gauss norm yields $\widehat{\mathbb{C}_{0,p}(z) \cap K[[z]]} = E_{0,K}$. Finally, it is clear that $\widehat{\mathbb{C}_{0,p}(z) \cap K[[z]]} = E_{0,p} \cap K[[z]]$.
\end{rema}

\begin{proof}[Proof of Proposition~\ref{prop_rational}]

Let $A$ be the companion matrix of $\mathcal{L}$, and assume that the order of $\mathcal{L}$ is equal to $n$. By assumption, we know that $\mathcal{L}$ has a strong Frobenius structure of period $h$. Since $K$ is a Frobenius field, the identity $i:K\rightarrow K$ map is a Frobenius endomorphism. Thus, there exists a matrix $H\in GL_n(E_{p})$ such that $\delta H=AH-p^hHA(z^{p^h})$. Since $\mathcal{L}$ is MOM at zero, it follows from Remark~\ref{rem_sol} that a fundamental matrix of solutions of $\delta X=AX$ is given by $Y_AX^{A(0)}$ with $Y_A\in GL_n(K[[z]])$ and $Y_A(0)=I$. Hence, a fundamental matrix of solutions of $\delta X=p^hA(z^{p^h})$ is given by $Y_A(z^{p^h})X^{p^hA(0)}$. Therefore, the equality $\delta H=AH-p^hHA(z^{p^h})$ implies that $HY_A(z^{p^h})X^{p^hA(0)}$ is a fundamental matrix of solution of $\delta X=AX$. Consequently, there exists $C\in GL_n(\mathbb{C}_p)$ such that $HY_A(z^{p^h})X^{p^hA(0)}=Y_AX^{A(0)}C$. Since $Y_A(0)=I$, $H\in M_n(E_{p})$, and $\log z$ is transcendental over $\mathbb{C}_p[[z,1/z]]$, we obtain $HY_A(z^{p^h})=Y_AC$. Hence, $H=Y_ACY_A(z^{p^h})^{-1}$. Since $Y_A\in GL_n(K[[z]])$, we conclude that  $H\in M_n(\bm{K}[[z]])\cap GL_n(E_p)$, where $\bm{K}$ is a finite extension of $K$ containing all the coefficients of $C$. In addition, we know that $E_p\subset\mathcal{A}_p$ and since $H\in M_n(\bm{K}[[z]])\cap GL_n(E_p)$, it follows from Proposition~1.2 of \cite{C86} that there exists $d\in \bm{K}$ such that $||H||=|d|$. Notice that $d\neq0$ because $H\in GL_n(E_{p})$. Let us set $T=\frac{1}{d}H$. Then, $T\in GL_n(E_{p})\cap M_n(\bm{K}[[z]])$, $||T||=1$, and $\delta T=AT-p^hTA(z^{p^h})$.

 Let $k\geq1$ be an integer. We set $T_k=T(z)T(z^{p^h})\cdots T(z^{p^{(k-1)h}})$. Thus, $T_k\in GL_n(E_p)\cap M_n(\bm{K}[[z]])$, $||T_k||\leq1$, and 
\begin{equation}\label{eq_k_h}
\delta T_{k}=AT_k-p^{kh}T_kA(z^{p^{kh}}).
\end{equation} 
 Since $f(z)$ is a solution of $\mathcal{L}$ and $A$ is the companion matrix of $\mathcal{L}$, the vector $(f,\delta f,\ldots, \delta^{n-1}f)^t$ is a solution of the system $\delta\vec{y}=A\vec{y}$. Thus, the vector  $(f(z^{p^{kh}}),(\delta f)(z^{p^{kh}}),\ldots, (\delta^{n-1}f)(z^{p^{kh}}))^t$ is a solution of the system $\delta\vec{y}=p^{kh}A(z^{p^{kh}})\vec{y}$.  Let us write $T_{k}=(t^{(k)}_{i,j})_{1\leq i,j\leq n}$.Then, thanks to Equation~\eqref{eq_k_h}, we deduce from Remark~\ref{rem_frob_ac} that $$t^{(k)}_{1,1}f(z^{p^{kh}})+t^{(k)}_{1,2}(\delta f)(z^{p^{kh}})+\cdots+t^{(k)}_{1,n}(\delta^{n-1}f)(z^{p^{kh}})$$
is a solution of $\mathcal{L}$. Notice that this solution belongs to $\bm{K}[[z]]$ because $T_k\in M_n(\bm{K}[[z]])$ and $f(z)\in 1+zK[[z]]$. Moreover,  this solution is different from zero because $(f(z^{p^h}),(\delta f)(z^{p^h}),\ldots, (\delta^{n-1}f)(z^{p^h}))^t$ is not the zero vector and $T_k$ belongs to $GL_n(E_{p})$. 

Since $\mathcal{L}$ is MOM at zero with coefficients in $E_{0,K}$, it follows from Remark~\ref{rem_sol} that there exists $c\in K$ such that 
\begin{equation*}
t^{(k)}_{1,1}f(z^{p^{kh}})+t^{(k)}_{1,2}(\delta f)(z^{p^{kh}})+\cdots+t^{(k)}_{1,n}(\delta^{n-1}f)(z^{p^{kh}})=cf.
\end{equation*}
By applying $\Lambda^{kh}_p$ to the previous equality, we have 
\begin{equation}\label{eq_lambda1}
\Lambda^{kh}_p(t^{(k)}_{1,1})f+\Lambda^{kh}_p(t^{(k)}_{1,2})\delta f+\cdots+\Lambda^{kh}_p(t^{(k)}_{1,n})\delta^{n-1}f=c\Lambda^{kh}_p(f).
\end{equation}
As, for all integers $j>0$, $(\delta^j f)(0)=0$ and $f(0)=1$, then $t^{(k)}_{1,1}(0)=c$. But, we know that $||T_k||\leq1$. Hence, $0<|c|\leq1$ and $c=\pi_{\bm{K}}^{r}u$, where $u$ is a unit of $\mathcal{O}_{\bm{K}}$, $\pi_{\bm{K}}$ is a uniformizer of $\mathcal{O}_{\bm{K}}$, and $r\in\mathbb{N}_{\geq0}$. Therefore, from Equation~\eqref{eq_lambda1}, we get 
\begin{equation}\label{eq_valor_propioI}
\frac{\Lambda^{kh}_p(t^{(k)}_{1,1})f+\Lambda^{kh}_p(t^{(k)}_{1,2})\delta f+\cdots+\Lambda^{kh}_p(t^{(k)}_{1,n})\delta^{n-1}f}{\pi_{\bm{K}}^r}=u\Lambda^{kh}_p(f).
\end{equation}
According to Corollary~\ref{prop_induction}, we know that $\frac{f(z)}{(\Lambda^{r+kh}_pf)(z^{p^{r+kh}})}$ belongs to $E_{0,K}$. Whence, there exists $P_{r+kh}(z)\in K_0(z)$ with norm 1 such that $P_{r+kh}[(\Lambda^{r+kh}_pf)(z^{p^{r+kh}})]\equiv f\bmod\pi^{r+kh}\mathcal{O}_{K}[[z]].$ Since, $\pi\mathcal{O}_K\subset\pi_{\bm{K}}\mathcal{O}_{\bm{K}}$, we also have
\begin{equation}\label{eq_kh}
P_{r+kh}[(\Lambda^{r+kh}_pf)(z^{p^{r+kh}})]\equiv f\bmod\pi_{\bm{K}}^{r+kh}\mathcal{O}_{\bm{K}}[[z]].
\end{equation}
Thus, for all integers $i>0$, $$(\delta^i P_{r+kh})[(\Lambda^{r+kh}_pf)(z^{p^{r+kh}})]\equiv\delta^{i}f\bmod\pi_{\bm{K}}^{r+kh}\mathcal{O}_{\bm{K}}[[z]].$$
In particular, for all integers $i\geq0$, there exists $d_i\in\mathcal{O}_K[[z]]$ such that $$(\delta^i P_{r+kh})[(\Lambda^{r+kh}_pf)(z^{p^{r+kh}})]-\delta^{i}f=\pi_{\bm{K}}^{r+kh}d_i\quad\text{ with }d_i\in\mathcal{O}_{\bm{K}}[[z]].$$
So, by using Equation~\eqref{eq_valor_propioI}, we obtain the following equalities
\begin{align*}
&\left[\frac{\sum_{i=0}^{n-1}\Lambda^{kh}_p(t^{(k)}_{1,i+1})\delta^iP_{r+kh}}{\pi_{\bm{K}}^r}\right](\Lambda^{r+kh}_pf)(z^{p^{r+kh}})-u\Lambda^{kh}_p(f)\\
&=\frac{\sum_{i=0}^{n-1}\Lambda^{kh}_p(t^{(k)}_{1,i+1})\left[\delta^{i}P_{r+kh}(\Lambda^{r+kh}_pf)(z^{p^{r+kh}})-\delta^{i}f\right]}{\pi_{\bm{K}}^r}\\
&=\frac{\sum_{i=0}^{n-1}\Lambda^{kh}_p(t^{(k)}_{1,i+1})(\pi_{\bm{K}}^{r+kh}d_i)}{\pi_{\bm{K}}^r}\\
&=\pi_{\bm{K}}^{kh}\left(\sum_{i=0}^{n-1}\Lambda^{kh}_p(t^{(k)}_{1,i+1})d_i\right).
\end{align*}
So
\begin{equation*}
\frac{\sum_{i=0}^{n-1}\Lambda^{kh}_p(t^{(k)}_{1,i+1})\delta^iP_{r+kh}}{\pi_{\bm{K}}^r}=\frac{\pi_{\bm{K}}^{kh}\left(\sum_{i=0}^{n-1}\Lambda^{kh}_p(t^{(k)}_{1,i+1})d_i\right)+u\Lambda^{kh}_p(f)}{(\Lambda^{r+kh}_pf)(z^{p^{r+kh}})}.
\end{equation*}
From (i) of Theorem~\ref{theo_analytic_element}, $f(z)\in 1+z\mathcal{O}_K[[z]]$. Thus, $(\Lambda^{r+kh}_pf)(z^{p^{r+kh}})$ is a unit of $\mathcal{O}_K[[z]]$ and from the previous equality we then have $$\frac{\sum_{i=0}^{n-1}\Lambda^{kh}_p(t^{(k)}_{1,i+1})\delta^iP_{r+kh}}{\pi_{\bm{K}}^r}\in\mathcal{O}_{\bm{K}}[[z]].$$
So, by reducing the last equality modulo $\pi_{\bm{K}}^{kh}$ and using the fact that $(\Lambda^{r+kh}_pf)(z^{p^{r+kh}})$ is a unit of $\mathcal{O}_K[[z]]$, we get that
\begin{equation}\label{eq_kh_mod}
\frac{\sum_{i=0}^{n-1}\Lambda^{kh}_p(t^{(k)}_{1,i+1})\delta^iP_{r+kh}}{\pi_{\bm{K}}^r}(\Lambda^{r+kh}_pf)(z^{p^{r+kh}})\equiv u\Lambda^{kh}_p(f)\bmod\pi^{kh}\mathcal{O}_{\bm{K}}[[z]].
\end{equation}
It follows from Equation~\eqref{eq_kh} that, $$P_{r+kh}[(\Lambda^{r+kh}_pf)(z^{p^{r+kh}})]\equiv f\bmod\pi_{\bm{K}}^{kh}\mathcal{O}_{\bm{K}}[[z]].$$
Since $f(0)=1=\Lambda^{r+kh}_p(f)(0)$, we obtain $P_{r+kh}(0)\equiv1\bmod\pi^{kh}$. Further, we know that $|P_{r+kh}|_{\mathcal{G}}=1$ and, by Remark~\ref{rem_k_0}, we obtain $P_{r+kh}\in\mathcal{O}_{K}[[z]]$ and thus, $P_{r+kh}$ is a unit of $\mathcal{O}_{K}[[z]]$. Whence, $P_{r+kh}\bmod\pi^{kh}\mathcal{O}_{K}[[z]]$  belongs to $1+z(\mathcal{O}_{K}/\pi^{kh})[[z]]$. So that, $P_{r+kh}\bmod\pi^{kh}\mathcal{O}_{K}[[z]]$ is a unit element of $(\mathcal{O}_{K}/\pi^{kh})[[z]]$. Consequently, $$(\Lambda^{r+kh}_pf)(z^{p^{r+kh}})=\frac{f}{P_{r+kh}}\bmod\pi^{kh}\mathcal{O}_{K}[[z]].$$
By substituting the above equality into Equation~\eqref{eq_kh_mod}, we get $$\frac{\sum_{i=0}^{n-1}\Lambda^{kh}_p(t^{(k)}_{1,i+1})\delta^iP_{r+kh}}{\pi_{\bm{K}}^r}\frac{f}{P_{r+kh}}\equiv u\Lambda^{kh}_p(f)\bmod\pi_{\bm{K}}^{kh}\mathcal{O}_{\bm{K}}[[z]].$$
Since $f(z)\in 1+z\mathcal{O}_K[[z]]$, $f\bmod\pi^{rh}\mathcal{O}_K[[z]]$ is a unit element of $(\mathcal{O}_K/\pi^{kh})[[z]]$. Then, from the previous equality, we have 
\begin{equation*}
\frac{\sum_{i=0}^{n-1}\Lambda^{kh}_p(t^{(k)}_{1,i+1})\delta^iP_{r+kh}}{\pi_{\bm{K}}^r}\frac{u^{-1}}{P_{r+kh}}\equiv\frac{\Lambda^{kh}_p(f)}{f}\bmod\pi_{\bm{K}}^{kh}\mathcal{O}_{\bm{K}}[[z]].
\end{equation*}

Therefore, \begin{equation}\label{eq_final}
\left|\frac{\sum_{i=0}^{n-1}\Lambda^{kh}_p(t^{(k)}_{1,i+1})\delta^iP_{r+kh}}{\pi_{\bm{K}}^r}\frac{u^{-1}}{P_{r+kh}}-\frac{\Lambda^{kh}_p(f)}{f}\right|_{\mathcal{G}}<|\pi_{\bm{K}}|^{kh}.
\end{equation}

%Since $(\pi^{kh})\mathcal{O}_K\subset(\pi^k)\mathcal{O}_K$, we obtain $$\frac{\sum_{i=0}^{n-1}\Lambda^{kh}_p(t^{(k)}_{1,i+1})\delta^iP_{r+kh}}{\pi^r}\frac{u^{-1}}{P_{r+kh}}\equiv\frac{\Lambda^{kh}_p(f)}{f}\bmod\pi^{k}\mathcal{O}_K[[z]].$$
%Now, we are going to see that there exists $Q_k\in\mathbb{C}_p(z)$ such that $$\frac{\sum_{i=0}^{n-1}\Lambda^{kh}_p(t^{(k)}_{1,i+1})\delta^iP_{r+kh}}{\pi^r}\frac{u^{-1}}{P_{r+kh}}\equiv Q_k\bmod\pi_{\bm{K}}^{kh}\mathcal{O}_{\bm{K}}[[z]].$$ We know that $P_{r+kh}$ has norm equal to 1 and that $P_{r+kh}(0)\equiv 1\bmod\pi^{kh}$. For this reason $P_{r+kh}(0)$ is a unit element of $\mathcal{O}_K$. But, we already know $P_{r+kh}$ belongs to $\mathcal{O}_K[[z]]$, and thus $P_{r+kh}$ is a unit element of $\mathcal{O}_K[[z]]$. So that, $1/P_{r+kh}\in K_{0}(z)$. 
As $t^{(k)}_{1,1}, t^{(k)}_{1,2},\ldots, t^{(k)}_{1,n}\in E_{p}\cap\bm{K}[[z]]$, by Proposition~5.1 of \cite{C86}, $\Lambda^{kh}_p(t^{(k)}_{1,1}), \Lambda^{kh}_p(t^{(k)}_{1,2}),\ldots, \Lambda^{kh}_p(t^{(k)}_{1,n})$ belong to $E_p\cap\bm{K}[[z]]$. Further, for all integers $i\geq0$, $\delta^{i}P_{r+kh}$ belongs to $K[[z]]$ because, we already know that $P_{r+kh}(z)$ is a unit of $\mathcal{O}_K[[z]]$. Therefore,  $$\frac{\sum_{i=0}^{n-1}\Lambda^{kh}_p(t^{(k)}_{1,i+1})\delta^iP_{r+kh}}{\pi_{\bm{K}}^r}\frac{u^{-1}}{P_{r+kh}}\in E_p\cap \bm{K}[[z]].$$ 
and, since $\left|\Lambda^{kh}_p(f)\big/f\right|_{\mathcal{G}}=1$, it follows from Equation~\eqref{eq_final} that $$\left|\frac{\sum_{i=0}^{n-1}\Lambda^{kh}_p(t^{(k)}_{1,i+1})\delta^iP_{r+kh}}{\pi_{\bm{K}}^r}\frac{u^{-1}}{P_{r+kh}}\right|_{\mathcal{G}}=1.$$
Therefore, $\frac{\sum_{i=0}^{n-1}\Lambda^{kh}_p(t^{(k)}_{1,i+1})\delta^iP_{r+kh}}{\pi_{\bm{K}}^r}\frac{u^{-1}}{P_{r+kh}}$ converges for all $x\in D_0$. So, $$\frac{\sum_{i=0}^{n-1}\Lambda^{kh}_p(t^{(k)}_{1,i+1})\delta^iP_{r+kh}}{\pi_{\bm{K}}^r}\frac{u^{-1}}{P_{r+kh}}\in E_{0,p}\cap\bm{K}[[z]].$$
Hence, according to Remark~\ref{rem_ast}, there exists $Q_{k}(z)\in\bm{K}_{0}(z)$ such that $$\left|\frac{\sum_{i=0}^{n-1}\Lambda^{kh}_p(t^{(k)}_{1,i+1})\delta^iP_{r+kh}}{\pi^r}\frac{u^{-1}}{P_{r+kh}}-Q_{k}(z)\right|_{\mathcal{G}}<|\pi_{\bm{K}}|^{kh}.$$ Thus, from Equation~\eqref{eq_final}, we get
\begin{equation}\label{eq_1_2_3}
\left|Q_{k}(z) -\frac{\Lambda^{kh}_p(f)}{f}\right|_{\mathcal{G}}<|\pi_{\bm{K}}|^{kh}.
\end{equation}
Finally, as $\left|\Lambda^{kh}_p(f)\big/f\right|_{\mathcal{G}}=1$ then, from \eqref{eq_1_2_3}, we deduce that $|Q_{k}(z)|_{\mathcal{G}}=1$ and $$Q_{k}(z)\equiv\frac{\Lambda^{kh}_p(f)}{f}\bmod\pi_{\bm{K}}^{kh}\mathcal{O}_{\bm{K}}[[z]].$$ 
That completes the proof.

\end{proof}

%By assuming Corollary~\ref{lemm_paso_base}, we can  prove (i) of Theorem~\ref{theo_analytic_element}.
\subsection{Proof of (ii) and (iii) of Theorem~\ref{theo_analytic_element}}

We are now ready to prove parts (ii) and (iii) of Theorem~\ref{theo_analytic_element}.

\begin{proof} Let $A$ be the companion matrix of $\mathcal{L}$. By hypotheses, we know that $\mathcal{L}$ has a strong  Frobenius structure and let $h$ be the period.

(ii) We first prove that $f(z)/f(z^{p^h})\in E_{0,K}$. Since $\mathcal{L}$ has a strong Frobenius structure, it follows from Propositions~4.1.2, 4.6.4, and 4.7.2 of \cite{C83} that there exists a matrix $U$ invertible with coefficients in the generic disc $D(t,1^{-})$ such that $\delta U=AU$. Now, the eigenvalues of $A(0)$ are all equal to zero because, by hypotheses, $\mathcal{L}$ is MOM at zero. It is clear that $A\in M_n(E_{0,K})$ because, by hypotheses, $\mathcal{L}$ is a differential operator with coefficients in $E_{0,K}$. Hence, $A$ belongs to $\mathcal{M}_K$. Further, by hypotheses, we know that the coefficients of $\mathcal{L}$ have Gauss norm less than or equal to 1. So, $||A||\leq1$. This puts us in a position to apply Corollary~\ref{prop_induction} and thus, for any integer $k\geq1$, $\frac{f(z)}{(\Lambda_p^{kh}(f))(z^{p^{kh}})}$ belongs to $E_{0,K}$. Consequently, for any integer  $k\geq1$, there is $P_{kh}(z)\in K_0(z)$ with norm 1 such that $$P_{kh}(z)(\Lambda_p^{kh}(f))(z^{p^{kh}})\equiv f(z)\bmod \mathfrak{\pi}^{kh}\mathcal{O}_K[[z]].$$
According to (i) of Theorem~\ref{theo_analytic_element}, $f(z)\in 1+z\mathcal{O}_K[[z]]$. Hence, we can reduce $f(z)$ modulo $\mathfrak{\pi}^{kh}$. So, $f(z)\bmod\mathfrak{\pi}^{kh}\in1+z\left(\mathcal{O}_K/\mathfrak{\pi}^{kh}\right)[[z]]$ is a unit element of the ring $\left(\mathcal{O}_K/\mathfrak{\pi}^{kh}\right)[[z]]$. Therefore, 
\begin{equation}\label{eq_pi}
P_{kh}(z)\frac{(\Lambda_p^{kh}(f))(z^{p^{kh}})}{f(z^{p^{kh}})}f(z^{p^{kh}})\equiv f(z)\bmod \mathfrak{\pi}^{kh}\mathcal{O}_K[[z]].
\end{equation}
But, according to Proposition \ref{prop_rational}, there exists a finite extension $\bm{K}$ of $K$ such that, for any integer $k\geq1$, there exists $Q_{k}(z)\in\bm{K}_{0}(z)$ with norm 1 such that 
\begin{equation}\label{eq_norma}
Q_{k}(z)\equiv\frac{\Lambda^{kh}_p(f)}{f}\bmod\pi_{\bm{K}}^{kh}\mathcal{O}_{\bm{K}}[[z]].
\end{equation}
where $\pi_{\bm{K}}$ is a uniformizer of $\mathcal{O}_{\bm{K}}$.

%We now prove that, for all integers $k\geq1$, there exists $B_k(z)\in\bm{K}_{0}(z)$ such that
%$$f(z)\equiv B_{k}(z)f(z^{p^{kh}})\bmod\mathfrak{\pi}_{\bm{K}}^{kh}\mathcal{O}_{\bm{K}}[[z]].$$

%Let $k\geq1$ be an integer, and let $\bm{L}_k$ be an infinite extension of $\bm{K}$ such that $\bm{L}_k$ containing all the coefficients of $Q_k(z)$. Let $\pi_{\bm{L}_k}$ be a uniformizer of $\bm{L}_k$. Then $\pi_{\bm{K}}\mathcal{O}_{\bm{K}}\subset\pi_{\bm{L}_k}\mathcal{O}_{\bm{L}_k}$, and thus, Equation~\eqref{eq_norma} implies that 
%\begin{equation}\label{eq_pi_l}
%Q_k(z)\equiv\frac{\Lambda_p^{kh}(f(z))}{f(z)}\bmod\pi_{\bm{L}_k}^{kh}\mathcal{O}_{\bm{L}_k}[[z]].
%\end{equation}
Since $\pi\mathcal{O}_{K}\subset\pi_{\bm{K}}\mathcal{O}_{\bm{K}}$, it follows from Equations~\eqref{eq_pi} and \eqref{eq_norma} that, for all integers $k\geq1$, \begin{equation}\label{eq_111}
f(z)\equiv B_{k}(z)f(z^{p^{kh}})\bmod\mathfrak{\pi}_{\bm{K}}^{kh}\mathcal{O}_{\bm{K}}[[z]].
\end{equation} 
With $B_{k}(z)=P_{kh}(z)Q_k(z^{p^{kh}})$. 
We have $B_{k}(z)\in\bm{K}_0(z)$ because $P_{kh}(z)\in K_0(z)$, $Q_{k}(z)\in\bm{K}_{0}(z)$, and $K\subset\bm{K}$. Further, $\left|B_k(z)\right|_{\mathcal{G}}=1$ because $\left|P_{kh}(z)\right|_{\mathcal{G}}=1=\left|Q_k(z)\right|_{\mathcal{G}}$.

We now proceed to prove that, for all integers $k\geq1$, $$\frac{f}{f(z^{p^h})}\equiv\frac{B_{(k+1)h}(z)}{B_{kh}(z^{p^h})}\bmod\mathfrak{\pi}_{\bm{K}}^{kh}\mathcal{O}_{\bm{K}}[[z]].$$

From \eqref{eq_111}, we have $f(z)\equiv B_{k}(z)f(z^{p^{kh}})\bmod\mathfrak{\pi}_{\bm{K}}^{kh}\mathcal{O}_{\bm{K}}[[z]]$. Hence,
\begin{equation}\label{eq_2}
f(z^{p^h})\equiv B_{k}(z^{p^h})f(z^{p^{(k+1)h}})\bmod\mathfrak{\pi}_{\bm{K}}^{kh}\mathcal{O}_{\bm{K}}[[z]].
\end{equation}
From \eqref{eq_111}, we know that
\begin{equation}\label{eq_3}
f(z)\equiv B_{k+1}(z)f(z^{p^{(k+1)h}})\bmod{\pi}_{\bm{K}}^{kh}\mathcal{O}_{\bm{K}}[[z]].
\end{equation}
Thus, from Equations \eqref{eq_2} and \eqref{eq_3}, we obtain $$\frac{f}{f(z^{p^h})}\equiv\frac{B_{k+1}(z)}{B_{k}(z^{p^h})}\bmod{\pi}_{\bm{K}}^{kh}\mathcal{O}_{\bm{K}}[[z]].$$

Consequently, we conclude that $f(z)/f(z^{p^h})=\lim\limits_{k\rightarrow\infty}\frac{B_{k+1}(z)}{B_{k}(z^{p^h})}$. So, $f(z)/f(z^{p^h})$ belongs to $E_p$. We now show that $f(z)/f(z^{p^h})\in E_{0,p}$. For this purpose we have to prove that, for all integers $k\geq1$, the rational function $\frac{B_{k+1}(z)}{B_{k}(z^{p^h})}$ belongs to $\mathbb{C}_{0,p}(z)$.  Let $k$ be a positive integer. It follows from Equation \eqref{eq_111} that $$f(z)- B_{k}(z)f(z^{p^{kh}})\in\pi_{\bm{K}}^{kh}\mathcal{O}_{\bm{K}}[[z]].$$ Thus, $1-B_{k}(0)\in\pi_{\bm{K}}^{kh}\mathcal{O}_{\bm{K}}$. Since the norm is non-Archimedean and $|\pi|<1$, $|B_{k}(0)|=1$. Further, we know that $|B_{k}(z)|_{\mathcal{G}}=1$ and thus, for all $x\in D_0$, $|B_{k}(x)|=|B_{kh}(0)|=1$. In particular, for all $x\in D_0$, $B_{k}(x)\neq0$. Therefore, $1/B_{k}(z)\in\bm{K}_0(z)$. Finally, as $B_{(k+1)}(z)$ belongs to $\bm{K}_0(z)$ then $\frac{B_{(k+1)}(z)}{B_{k}(z^{p^h})}$ belongs to $\bm{K}_0(z)$.  Consequently, $f(z)/f(z^{p^h})\in E_{0,p}$. 

But, it is clear that $f(z)/f(z^{p^h})\in K[[z]]$. So, $f(z)/f(z^{p^h})\in E_{0,p}\cap K[[z]]$ and thus, according to Remark~\ref{rem_ast}, we conclude that $f(z)/f(z^{p^h})\in E_{0,K}$.
\medskip

(iii) We now prove that $f'(z)/f(z)\in E_{0,K}$. We put $H(z)=f(z)/f(z^{p^h})$.  According to (ii), $H$ belongs to $E_{0,K}$. Let us show that $H$ is a unit of $E_{0,K}$. For this purpose, we need to show that, for all $x\in D_0$, $H(x)\neq0$. We know from (i) of Theorem~\ref{theo_analytic_element} that $f(z)\in 1+z\mathcal{O}_K[[z]]$. Since the norm is non-Archimedean, we get that, for all $x\in D_0$, $|f(x)|=1$. Whence, for all $x\in D_0$, $|H(x)|=1$ and therefore, $H(x)\neq0$ for all $x\in D_0$. Now, let us show that $H'(z)$ also belongs to $E_{0,K}$. Since $H(z)\in E_K$, $H'(z)\in E_K$. So, it is sufficient to show that, for all $x\in D_0$, $H'(x)$ converges. It is not hard to see that $$H'(z)=H(z)\left[\frac{f'(z)}{f(z)}-p^{h}z^{p^{h}-1}\frac{f'(z^{p^{h}})}{f(z^{p^{h}})}\right].$$ 
As $f(z)\in 1+z\mathcal{O}_K[[z]]$ then, for all $x\in D_0$, $f(x)$ converges. So, from the previous equality, we conclude that, for all $x\in D_0$, $H'(x)$ converges. For every $m\geq0$, we consider $H_m=H(z)H(z^{p^h})\cdots H(z^{p^{mh}})$. So, $H_m\in E_{0,K}$ and is a unit of $E_{0,K}$. Further, we also have
\begin{equation}\label{eq_frob_m}
H'_m(z)=H_m(z)\left[\frac{f'(z)}{f(z)}-p^{(m+1)h}z^{p^{(m+1)h}-1}\frac{f'(z^{p^{(m+1)h}})}{f(z^{p^{(m+1)h}})}\right].
\end{equation}
Since $f(z)\in 1+z\mathcal{O}_K[[z]]$, we get that, for every integer $m\geq0$, $\frac{f'(z^{p^{(m+1)h}})}{f(z^{p^{(m+1)h})}}\in\mathcal{O}_K[[z]]$ and $H_m\in 1+z\mathcal{O}_K[[z]]$. As $p\in (\pi)\mathcal{O}_K$ then, it is derived from Equation~\eqref{eq_frob_m} that, for every integer $m\geq0$, $$\frac{H'_m(z)}{H_m(z)}\equiv\frac{f'(z)}{f(z)}\bmod\pi^{(m+1)h}\mathcal{O}_K[[z]].$$
But, for every integer $m\geq0$, $\frac{H'_m(z)}{H_m(z)}\in E_{0,K}$. The ring $E_{0,K}$ is complete because $K$ is complete and thus, we have $f'(z)/f(z)\in E_{0,K}$. This completes the proof.

\end{proof}

\end{document}